%% file: coasch1.tex
\documentclass[reqno]{amsart}
\usepackage[top=2cm,bottom=2cm,right=2cm,left=2cm]{geometry}

\input{preamble.tex}

\title[Empirical measures for maps with $1$-regularly varying return time distribution]{Almost sure behaviour of empirical measures for intermittent maps whose return time distribution is $1$-regularly varying}

\author{Douglas Coates and Tanja I.~Schindler}

\begin{document}

\maketitle

\begin{abstract}
  In this article we study sets of limit points $\cL$ for the sequence of empirical measures of intermittent maps with several equally sticky neutral fixed points. Such
  intermittent maps come with a parameter $\alpha > 0$ which determines the ``stickiness'' of the neutral fixed points.
  When $\alpha \neq 1$, it is known that two different behaviours can occur: either $\cL$ is almost surely the singleton $\{\mu\}$ where $\mu$ is an absolutely continuous invariant probability measure; or $\cL$ is almost surely given by the entire simplex $\pi\cS$ of invariant measures supported at the neutral fixed points.
  We study the case that $\alpha = 1$ and show that neither of these two behaviours occur and that $\cL$ is almost surely equal to a star-shaped subset of the simplex $\pi \cS$. Moreover, we show that when slowly varying functions are introduced one can obtain a large range of different behaviours, including both the possible outcomes when $\alpha \neq 1$.
  Our results apply to a large class of intermittent interval maps as well as parabolic rational maps of the Riemann sphere. 
\end{abstract}


\section{Introduction}
\label{sec:intro}

In this article we set out to understand the typical behaviour of the sequences of empirical measures for a class of intermittent dynamical systems. Intermittent transformations were first introduced in \cite{PomMan1980} in order to model the transition between laminar and turbulent flows in fluids. Mathematically, they are maps that are expanding everywhere except at some finite number of neutral fixed points. Their dynamics exhibit \emph{intermittent} behaviour that is characterised by short chaotic bursts (where the orbit stays in the expanding region of phase space) separated by long predictable laminar phases (where the orbit stays close to a neutral fixed point). These are a widely studied and important class of dynamical systems as they provide one of the simplest settings to study non-uniformly hyperbolic phenomena and are often viewed as a test-bed for more complex systems. There is a plethora of literature on these systems (see e.g.\ \cite{LivSauVai1999,Sar2002,Gou2004,Gou2004a,MelTer2012,Gou2011,BalTod2016,Kor2016} and the references therein).

An intermittent transformation $f \colon X \to X$ of a probability space $(X,m)$ comes with a natural parameter $\alpha > 0$ which describes how intermittent the system is, or equivalently how ``sticky'' the neutral fixed points are. Under some quite mild hypothesis (see e.g.\ \cite{Zwe1998,BandtJen}) an intermittent system will preserve an
ergodic measure $\mu$
that is equivalent to $m$.
If the neutral fixed points are not too sticky, which will correspond to $\alpha > 1$, this measure can be normalised to a probability measure. If instead $\alpha \in (0, 1]$ then $\mu$ is only $\sigma$-finite. In this case the ergodic theorem for $\mu$ becomes degenerate in the sense that $\lim_{ n \to \infty} \frac{ 1 }{ n } \sum_{ k = 0 }^{ n - 1 } u \circ f^k \aseq 0$ for every $u \in L^1 ( \mu )$
where by ``$\aseq$'' we denote almost sure equality with respect to $m$.\footnote{
  As the measures $m$ and $\mu$ are equivalent $a \aseq b$ is the same as $a = b$ $\mu$ almost everywhere.
}
exist. Let
\begin{equation*}
  e_n \colon X \to \cM_1 (X),
  \quad
  e_n (x) \coloneqq \frac{1}{n} \sum_{ k = 0 }^{ n - 1 } \delta_{ f^k x }
\end{equation*}
denote the sequence of \emph{empirical measures} of $f$. Notice that if $e_n(x)$ converges to some $\nu$ for a given $x \in X$, then, by definition, the limit describes the asymptotic distribution of $x$ in phase space.
We will say that a probability measure $\nu$ is \emph{physical} if it describes the asymptotic distribution of a large set of initial conditions in the sense that $m \left\{ x \colon e_n (x) \to \nu \right\} > 0$.
The study of empirical measures is a classical and central topic in ergodic theory.
The existence and non-existence of physical measures and more generally the behaviour of empirical measures for intermittent maps has been studied recently in \cite{AarThaZwe2005,CoaLuz2024,CoaMelTal2026}. Presently there is a reasonably complete understanding of the almost sure behaviour of the sequence $e_n$ in the case that $\alpha \neq 1$. When $\alpha \neq 1$ one can give a precise description of the set of limit points
\begin{equation*}
  \cL ( x ) \coloneqq \{ \text{limit points of } e_n (x) \},
\end{equation*}
for almost every initial condition $x \in X$.
In this article we aim to complete this picture by treating the case $\alpha = 1$. In order to continue this discussion and to present our main results it is convenient to fix a specific class of intermittent systems. Here in the introduction we will state our results in the context of interval maps and parabolic rational maps but we note that these results are consequences of more abstract results that we present later in Section \ref{sec:results}.

\subsection{Thaler's intermittent maps of the interval}

Let $X = [0,1]$ and let $m = \Leb$. We let $f \colon X \to X$ be a full-branched map with finitely many neutral fixed points and some number (possibly countably many) of uniformly expanding branches. Specifically, we suppose that:
\begin{description}
  \item[T1\label{itm:T1}] $\cP$ is a partition of $X$ (modulo $m$) into open sub-intervals and the restriction $f_{I} = f|_{I} \colon I \to (0,1)$ is a twice differentiable bijection for each $I \in \cP$;

  \item[T2] there exist $d \ge 1$ neutral fixed points $\xi_1, \ldots, \xi_d$ with $f'\xi_i = 1$ for all $i = 1,\ldots, d$ and $f' x > 1$ for every $x\not\in \{ \xi_1, \ldots, \xi_d \}$;

  \item[T3\label{itm:T3}] the map satisfies \emph{Adler's condition}: $\displaystyle \sup_{ I\in \cP } \| f''_I / (f'_I)^2 \|_{\infty} < \infty$;

 \item[T4\label{itm:expansion}] there exist an $\alpha \in (0,1]$ and $a_1, \ldots, a_d>0$ such that for every $i = 1, \ldots, d$,  $| f x - \xi_i | \sim a_i | x - \xi_i |^{ 1 + 1/\alpha}$.
\end{description}

\textbf{Notation.} Throughout we will use the $x_{ n } \sim y_{ n }$ to denote $x_{ n } = y_{ n } ( 1 + o(1))$ and $x_{ n } \ll y_{ n }$ if $x_{ n } = O(y_{ n })$. We will also write $x_{ n } \approx y_{ n }$ if both $y_{ n } \ll x_{ n }$ and $x_n \ll y_{ n }$.

Such maps were introduced in \cite{Tha1980,Tha1983} and were shown to be conservative, ergodic and preserve a $\sigma$-finite ergodic measure $\mu$ which is equivalent to Lebesgue and satisfies $\mu(X) = \infty$. 
If $\alpha \in (0, 1)$, then it was shown in~\cite{AarThaZwe2005,CoaLuz2024} in the case that $d = 2$ and then in~\cite{CoaMelTal2026} for general $d$, that
the set of limit points $\cL$
is almost surely equal to the whole simplex of invariant measures
supported at the neutral fixed points:
\begin{equation}\label{eq:as-lim-alpha-not-1}
  \cL \aseq \pi \cS,
  \quwhere
  \cS\coloneqq \{ \vp \colon p_1 + \dots + p_d = 1,\;p_i  \geq 0 \},
  \quand
  \pi \colon \cS \to \cM_1 ( X ),
  \quad
  \pi ( \vp ) = \sum_{ i = 1 }^{ d } p_{ i } \delta_{ \xi_i }.
\end{equation}
When $\alpha = 1$, the methods in~\cite{AarThaZwe2005,CoaLuz2024,CoaMelTal2026} yield that there exists a $\vrho \in \cS$ so that
\[
  \{ \delta_{\xi_1} , \ldots, \delta_{ \xi_d  }, \pi(\vrho)\}\subset \cL
  \subset
  \pi \cS,
\]
almost surely. Moreover it is not hard to see that the set of limit points is necessarily connected. This means that in the case that $d \le 2$ one again obtains that $\cL \aseq \pi\cS$. We will show that this is no longer the case when $d \ge 3$. Our main theorem, in the context of the present example, states:

\begin{theorem}\label{thm:theorem-constant}
  Suppose that $f$ satisfies \ref{itm:T1}--\ref{itm:expansion} with $\alpha = 1$. Then, there exists a $\vrho \in \cS$ so that $\cL \aseq  \pi \left(\bigcup_{ i=1 }^{ d } \RR_i\right)$, where $\RR_i$ is the straight line that connects the point $\vrho$ to the point $(0,\ldots,0,1,0,\ldots,0)$, the $i^{ th }$ vertex of the simplex. 
\end{theorem}

\begin{remark}
  In Figure \ref{fig:1trim} we show a depiction of the subset $\bigcup_{ i=1 }^{ d } \RR_i$ in the case that $d = 4$.
\end{remark}

\begin{remark}
  For any $\alpha \in (0, 1]$  it was shown in \cite{CoaMelTal2026} that despite the fact that the $e_n$ do not converge almost surely, they do converge in distribution. This captures exactly how the $e_n$ are asymptotically distributed in $\pi\cS$. When $\alpha \in (0,1)$ this limiting distribution is non-degenerate and gives positive mass $\pi U$ for any non-empty open subset $U \subset \cS$. The non-degeneracy of this limiting distribution is a key step in the approach in \cite{CoaMelTal2026} for proving \eqref{eq:as-lim-alpha-not-1}.

  When $\alpha = 1$ the limiting distribution of the $e_n$ is degenerate, in fact it is the point mass concentrated at a point $\pi ( \vrho )$. In this case the approach in \cite{CoaMelTal2026} for establishing $\cL$ fails completely. In this article we develop a different approach which hinges on the trimmed strong law of large numbers proved in \cite{AarNak2003}.
\end{remark}

\begin{remark}\label{rem:infinite}
  The properties \ref{itm:T1}--\ref{itm:T3} already imply that the invariant measure is infinite and so it does not make sense to talk about the $\alpha > 1$ case in this context. Indeed,  \ref{itm:T1}--\ref{itm:T3} combined with Taylor's Theorem ensure that \ref{itm:expansion} cannot hold with $\alpha > 1$.
  However, if instead of \ref{itm:T3} and \ref{itm:expansion} we had assumed that $f''(x) < \infty$ for every $x \not\in \{ \xi_1, \ldots, \xi_d \}$ and
  \[
    | f'' (x) | \sim \frac{ 1 }{ \alpha } \left( 1 + \frac{ 1 }{ \alpha } \right) a_i | x - \xi_i |^{ 1 / \alpha - 1 },
    \text{ as } x \to \xi_i
  \]
  for every $i = 1, \ldots, d$ for some $\alpha > 0$ and some $c_1, \ldots, c_d$ then it would make perfect sense to consider $\alpha > 1$. For $\alpha \in (0,1]$ this case would just be a special case of the one described above. For $\alpha > 1$ we  would instead find the measure $\mu$ would be finite and that $\cL \aseq \{\mu\}$.
\end{remark}

In the next section we generalize the above setting by introducing additionally a slowly varying function in \eqref{eq:as-lim-alpha-not-1}.

\subsection{Thaler maps with slowly varying functions}
\label{sec:thaler-slow}

When $\alpha = 1$ and the constants $a_1,\ldots,a_d > 0$ in \ref{itm:expansion} are replaced by functions $\beta_1,\ldots, \beta_d \colon [0,1] \to \R$ which are slowly varying many different subsets of $\pi \cS$ can appear as the set of almost sure limit points of the empirical measures. Recall that a function $\ell \colon \R \to \R$ is \emph{slowly varying} (at a point $p \in \bar \R$) if $\lim_{ x \to p } \ell ( \lambda x ) / \ell ( x ) = 1$ for every $\lambda > 0$. A classical example is $\log$ which is slowly varying at $\infty$.
To that end, let us now assume that \ref{itm:expansion} is replaced by the following:
\begin{description}
  \item[T4$^\prime$\label{itm:exp-slow}] there exist $\beta_1,\ldots, \beta_d \colon [0,1] \to \R$ that are real-analytic functions and slowly varying at $0$ such that for every $i = 1, \ldots, d$,  $| f x - \xi_i | \sim \beta_i ( | x - \xi_i| )  | x - \xi_i |^{ 2 }$.
\end{description}
In order to state our results in this more general setting we need to introduce some additional notation and terminology.

In order to determine the set $\cL$ it is convenient to introduce for each $i = 1, \ldots, d$ the function
\begin{equation}\label{eq:def-elli-intro}
  \ell_i \colon \R^+ \to \R^+,
  \quad
  \ell_i ( x ) = c_i \beta_i^\sharp ( 1 / x )
\end{equation}
where $c_i > 0$ is a positive constant and where $\beta_i^\sharp$ is the \emph{de Bruijn conjugate}\footnote{
  For a slowly varying function $\ell$, its de Bruijn conjugate $\ell^\sharp$ is the slowly varying function such that
  \[
    \lim_{ x \to p } \ell ( x ) \ell^\sharp ( x  \ell ( x ) ) = 1,
    \quand
    \lim_{ x \to p } \ell^\sharp ( x ) \ell ( x  \ell^\sharp ( x ) ) = 1.
  \]
  The de Bruijn conjugate always exists and is unique up to asymptotic equivalence (see e.g.\ \cite[Theorem 1.5.13]{BinGolTeu1987}). Notice that hence the function $\ell_i$ defined in \eqref{eq:def-elli-intro} is slowly varying at $\infty$.
} of $\beta_i$. The precise value of $c_i$ is given later in \eqref{eq:def-ell-tha}.
We need to impose the following (reasonably mild) condition on the $\ell_i$, namely that the limit
\begin{equation}\label{eq:rho_i}
  \rho_i = \lim_{ x \to \infty } \frac{ \ell_i (x) }{ \ell_1 ( x ) + \cdots + \ell_d ( x ) } \in [0,1], \quad \text{exists for all} \quad i = 1, \ldots, d.
\end{equation}
Notice that whether or not this limit exists is determined only by the functions $\beta_i$ and not by the value of the constants $c_i$. Moreover, if $\rho_i \in \{0,1\}$ then changing the value of $c_i$ would not change the value of the corresponding $\rho_i$.

\begin{remark}
  In the special case that each of the functions $\beta_i$ are constant their de Bruijn conjugates are also constant, and so it is clear that the $\rho_i$ all exist. In this case the $\rho_i$ are precisely the components of the vector $\vrho$ which appears in Theorem \ref{thm:theorem-constant}.
\end{remark}

In light of the previous remark, we define $\vrho \in \cS$ to be the vector whose components are given by $\rho_i$. As before we let $\RR_i$ be the straight line, contained in $\cS$ which joins the $i^{ th }$ vertex of the simplex with the point $\vrho$. Moreover, for $N \ge 1$ and $\zeta \in \{1, \ldots, d \}^{ N }$ we let $\RR_{ \vz }$ be the convex hull of the lines $\RR_{\zeta_1}, \ldots, \RR_{\zeta_N}$.

We also need to define
\begin{equation}\label{eq:def-Li-Ni}
  L_i (x) \coloneqq \int_{ 0 }^x  \frac{ \ell_i ( t ) }{ t }  \, dt,
  \quand
  \mathfrak{N}^{(i)} \coloneqq \min \left\{ 
    k \ge 1 \colon \sum_{ n = 1 }^{ \infty }
    \left(\frac{ n \ell_i (n) }{ L_i (n) } \right)^{ k + 1 } \frac{ 1 }{ n } < \infty
  \right\}.
\end{equation}
The quantities $L_i,\NN^{(i)}$ give additional measures of how ``sticky'' the corresponding fixed point is.
It turns out that if $L_i$ is bounded then the corresponding fixed point is not very sticky and will play no role in determining the limit set $\cL$. For this reason we will make the simplifying assumption\footnote{The assumptions \ref{itm:T1}--\ref{itm:T3} ensure that the invariant measure $\mu$ is infinite (see Remark \ref{rem:infinite}). As will become clear in Section \ref{subsec: inter-inter-maps}, this in turn implies that at least one the $L_i$ must be undbounded.} that
\begin{equation}\label{eq:ellnonint}
  \lim_{ x \to \infty } L_i ( x ) = \infty.  
\end{equation}
Under this assumption the value of the $\NN_i$ allows us to determine the relative strengths of the fixed points. Roughly speaking,
the slower the function $\beta_i (| x - \xi_i |) | x - \xi_i |$ converges to zero, the larger the value of $\NN^{(i)}$ is and the stickier the fixed point $\xi_i$ is. 

\begin{theorem}\label{thm:N}
  Suppose that $f$ satisfies \ref{itm:T1}--\ref{itm:T3} and \ref{itm:exp-slow} and that \eqref{eq:rho_i} and \eqref{eq:ellnonint} hold. If $N = \NN^{(1)} = \ldots = \NN^{(d)}$ for some $1 \le N < \infty$, and $\vrho \in \text{int} \cS$ then 
    \[
      \cL \aseq \pi \left(\bigcup_{ \vi \in \{ 1, \ldots,  d \}^N } \RR_{\vi}\right).
    \]
\end{theorem}

\begin{figure}\label{fig}
  \centering
  \begin{subfigure}[t]{0.15\textwidth}
      \centering
    \includegraphics[width=\textwidth]{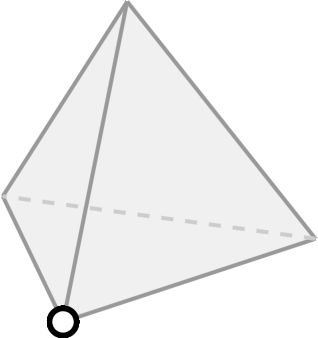}
    \caption{$\ell_1 = \log$, $\ell_2 = \cdots = \ell_4 = 1$.}
    \label{fig:phys}
  \end{subfigure}
  \hspace{2em}
  \begin{subfigure}[t]{0.15\textwidth}
      \centering
    \includegraphics[width=\textwidth]{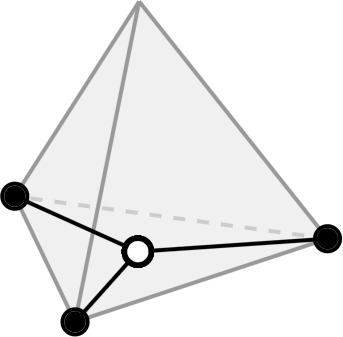}
    \caption{$\ell_1 = \cdots = \ell_3 = \log$, $\ell_4 = 1$.}
    \label{fig:weird}
  \end{subfigure}
  \hspace{2em}
  \begin{subfigure}[t]{0.15\textwidth}
      \centering
    \includegraphics[width=\textwidth]{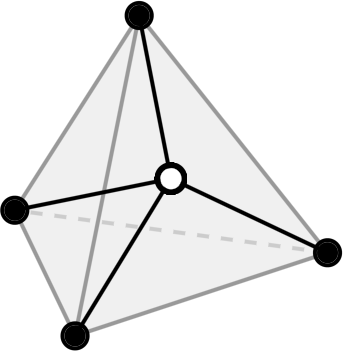}
    \caption{$\NN^{(1)} = \cdots = \NN^{(4)} = 1$ and $\vrho \in \operatorname{int} \cS$}
    \label{fig:1trim}
  \end{subfigure}
  \hspace{2em}
  \begin{subfigure}[t]{0.15\textwidth}
      \centering
    \includegraphics[width=\textwidth]{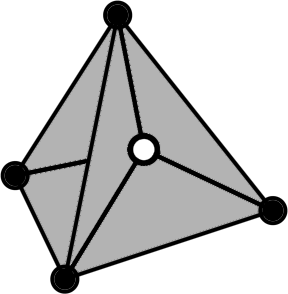}
    \caption{$\NN^{(1)} = \cdots = \NN^{(4)} = 2$ and $\vrho \in \operatorname{int} \cS$}
    \label{fig:2trim}
  \end{subfigure}
  \hspace{2em}
  \begin{subfigure}[t]{0.15\textwidth}
      \centering
    \includegraphics[width=\textwidth]{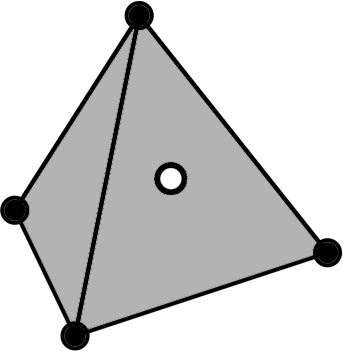}
    \caption{$\NN^{(1)} = \cdots = \NN^{(4)} = 3$ and $\vrho \in \operatorname{int} \cS$}
    \label{fig:3trim}
  \end{subfigure}
  \caption{Depictions of the set $\cL$ for a map satisfying \ref{itm:T1}--\ref{itm:T3} and \ref{itm:exp-slow} with $d = 4$. The figures show how $\cL$ changes for some different choices of $\beta_1,\ldots,\beta_4$ in \ref{itm:exp-slow}. Each figure shows the simplex $\cS$ in light grey, and the set $\pi^{-1} \cL$ with its interior coloured dark grey and its boundary coloured black. The point $\vrho$ is shown as a black circle with a white centre.}
\end{figure}

\begin{remark}
  \begin{itemize}
    \item
    If all of the $\NN^{(i)}$ are equal to $1$, then one finds the same behaviour as in Theorem \ref{thm:theorem-constant}. Conversely, if each of the $\beta_i$ are constant, then one finds that  \eqref{eq:rho_i} holds with each $\rho_i > 0$ and that all of the $\NN^{(i)} = 1$ (see Example~\ref{ex: ellsimc}). In particular, Theorem \ref{thm:theorem-constant} is just a special case of Theorem \ref{thm:N}.

    \item 
    If all of the $\NN^{(i)}$ are equal to some $2 \le N < d-1$ then, one obtains that $\cL$ is almost surely equal to a proper subset of $\pi\cS$ that differs from all cases previously considered, see Figure \ref{fig:2trim} for a sketch in the case that $d = 4$ and each $\NN^{(i)} = 2$.
    
    \item If, on the other hand, all of the $\NN^{(i)}$ are equal to some $N \geq d-1$ then, one obtains that $\cL$ is almost surely equal to the full set $\pi\cS$, see Figure \ref{fig:3trim} for a sketch in the case that $d = 4$ and each $\NN^{(i)} = 3$.
    
    \item
    The conditions of Theorem \ref{thm:phys-meas-case} are non-empty, in particular
    for any $N \ge 1$ one can always choose the $\beta_i$ in such a way that $N = \NN^{(i)}$  and $\rho_i > 0$.
  \end{itemize}
\end{remark}

We conclude this Section by giving an example of what different behaviours can occur when
some of the $\rho_i$ are zero.

\begin{theorem}
  \label{thm:phys-meas-case}
  Suppose that $f$ satisfies \ref{itm:T1}--\ref{itm:T3} and \ref{itm:exp-slow}, and that
  $\beta_i ( x ) = c_i /  | \log x  |$ for $i = 1, \ldots, j$ and that
  $\beta_i \equiv c_i$ for $i = j + 1, \ldots, d$ for some $c_1, \ldots, c_d > 0$.
  Then,
  \begin{equation*}
    \cL \aseq \pi \left(\bigcup_{ i=1}^j \RR_i\right).
  \end{equation*}
\end{theorem}

\begin{remark}
  Under the hypothesis of Theorem \ref{thm:phys-meas-case}, one finds that
  \eqref{eq:rho_i} holds and that each $\NN^{(i)} = 1$. However, Theorem \ref{thm:N} does not apply (unless $j = d$) as  $\rho_i > 0$ only for $i = 1, \ldots, j$, and hence $\vrho \in \partial \cS$. A special case of this theorem worth noting is that if $j = 1$ then there is a unique physical measure which is the Dirac mass at fixed point $\xi_1$. See Figures \ref{fig:phys} and \ref{fig:weird}.
\end{remark}

\begin{remark}
  Theorems \ref{thm:N} and \ref{thm:phys-meas-case} are special cases of Theorems \ref{thm:2-trimming} and \ref{thm: phys meas} respectively, which we will state in Section \ref{sec:results}.
\end{remark}

\subsection{Parabolic Rational Maps}\label{sec:parabolic}
We now turn our attention to parabolic rational maps of the Riemann sphere $\overline \C$. These are rational maps whose Julia set contains  parabolic periodic orbits but does not contain critical points.
We will assume:

\begin{description}
  \item[P1\label{p1}] $f \colon X \to X$ is a parabolic rational map acting on its Julia set $X \subset \overline\C$ which contains $d \ge 1$ indifferent periodic orbits.
\end{description}

Given a function $f$ satisfying \ref{p1} we will let $\dim X$ denote the Hausdorff dimension of $X$, and let $m$ be the unique
$\dim X$-conformal probability measure on $X$.
It was shown in \cite{AarDenUrb1993} that $f$ is a conservative and ergodic transformation of $(X,\mu)$ and that $f$ preserves a unique (up to scaling) $\sigma$-finite measure $\mu$ which is equivalent to the conformal measure $m$.

In what follows, for $i = 1, \ldots, d$, we will write $\nu_i$ for the invariant measure supported on the $i^{th}$ periodic orbit. Let $\cS$ be the simplex given in \eqref{eq:as-lim-alpha-not-1}, and define 
\[
  \Pi \colon \cS \to \cM_1 ( X ),
  \quad
  \Pi ( \vp ) = p_1 \nu_1 + \cdots + p_d \nu_d
\]
to be the natural identification between the simplex $\cS$ and the space of convex combinations of the measures $\nu_1, \ldots, \nu_d$.

It will often be convenient to reduce the periodic orbits to fixed points by passing to an iterate of $f$. To that end, let $M$ be a common period of the indifferent periodic orbits, let $g \coloneqq f^M $ and let $\xi_1, \ldots, \xi_{ d' }$ denote the indifferent fixed points of $g$ where $d' \geq d$. 
As $g$ is holomorphic we know that there exist $a_i, \gamma_i > 0$ such that
\begin{equation}\label{eq:g}
  g (z) = z + a_i  ( z - \xi_i )^{ 1 + \gamma_i } + O ( (z - \xi_i )^{ 2 + \gamma_i } )
\end{equation}
as $z \to \xi_i$ for each $i = 1, \ldots, d'$.
It is also shown in \cite{AarDenUrb1993} that whether or not the measure is finite depends on the relative value of $\gamma = \max \{ \gamma_1, \ldots, \gamma_{d'} \}$ and the Hausdorff dimension $\dim X$ of the Julia set.
By \cite[Corollary 8.6 and Theorem 8.8]{AarDenUrb1993}, $\dim X \in\left( \frac{ \gamma }{ \gamma + 1 }, 2\right)$ and by \cite[Theorem 9.8]{AarDenUrb1993}
\[
  \mu ( X ) = \infty
  \Longleftrightarrow
  \dim X \in
  \left( \frac{ \gamma }{ \gamma + 1 }, \frac{ 2 \gamma }{ \gamma + 1 }
  \right].
\]
We set $\alpha = ( 1 + 1 / \gamma ) \dim X - 1$, so that $\mu( X ) = \infty$ if and only if $\alpha \in (0,1]$. As only the parabolic fixed points of $g$ of maximal exponent will play a role, we will make the simplifying assumption that all of the exponents are equal. Moreover, as we are interested only in the case that $\alpha = 1$ we will also impose:
\begin{description}
  \item[P2\label{p2}] The powers $\gamma_i$ in \eqref{eq:g} satisfy $\gamma_1 = \cdots = \gamma_{d'} = \gamma$, where $\gamma = \left( \frac{ 2 }{ \dim X } - 1\right)^{ -1 }$ (so that $\alpha = ( 1 + 1 / \gamma ) \dim X - 1 = 1$). 
\end{description}

\begin{remark}
  One can check that \ref{p2} does not depend on the choice of $M$ in the definition of $g$.
\end{remark}

\begin{theorem}\label{thm:parabolic}
  Suppose that $f$ satisfies \ref{p1} and \ref{p2}. Then, there exists a $\vrho \in \operatorname{int}\cS$ (whose value is given in \eqref{eq:def-rho-para} below) such that
  \[
    \{ \text{limit points of } e_n \}
    \aseq
    \Pi \left( \bigcup_{ i = 1 }^{ d } \RR_i \right).
  \]
\end{theorem}

\begin{remark}
  Similarly to the case of intermittent maps we know that $\cL \aseq \{ \mu \}$ when $\alpha > 1$ (and hence $\mu$ is finite) and $\cL \aseq \Pi ( \cS )$ when $\alpha \in (0,1)$. The former is a simple consequence of the ergodic theorem and the latter can be deduced from \cite{AarDenUrb1993,CoaMelTal2026}.
\end{remark}

\subsection{An outline of the paper}

We begin in Section~\ref{sec:results} by giving statements of our main results in a more abstract setting. Then, in Section~\ref{sec:cf mixing} we introduce the concept of continued fraction mixing and trimming. 
In Section~\ref{sec:proofs}  we proceed to give the proofs of the results stated in Section~\ref{sec:results} and finally, we give the proofs for the statements regarding interval maps and parabolic rational maps in Sections~\ref{subsec: inter-inter-maps}
and~\ref{subsec:para-rat-map}.

\section{Assumptions and statement of the abstract results}
\label{sec:results}

Throughout we will assume that $(X, \cB, m)$ is a standard probability space. The abstract setting for which we prove our main results may be stated as follows.

\begin{description}
  \item[A1\label{itm:f-cons}]
  \textbf{Infinite measure preserving transformation:}
    $f: X \to X$  conservative transformation of a metric space $X$ which preserves a $\sigma$-finite ergodic measure  $\mu$ that is equivalent to $m$ and that satisfies $\mu( X ) = \infty$.

    \item[A2\label{itm:fxd-pnts}]\textbf{Measure is concentrated at fixed points:}
      There exist $d \ge 1$ fixed points $\xi_1, \ldots, \xi_d$ such that
      \[
        \mu \big[ X \setminus ( B_{\eps} (\xi_1)  \cup \cdots \cup B_{\eps} (\xi_d) ) \big] < +\infty,
        \quad \text{for all} \quad \eps > 0.
      \]

      \item[A3\label{itm:dyn-sep}]\textbf{Dynamical separation:}
      There exists a  subset $Y \subset X$ with $\mu (Y) = 1$
      which \emph{dynamically separates} the fixed points $\xi_1, \ldots, \xi_d$. In particular, there exists a partition (modulo $m$) of $X' \coloneqq X \setminus Y$ into open sets  $X_1 \cup \cdots \cup X_d$ such that for all $i = 1, \ldots, d$, $\xi_i \in X_i$, and for all $k \ge 1$, and $j \neq i$,
      \begin{equation}
        \label{eq:dyn-sep}
        x \in X_i \text{ and } f^k x \in X_j  \Longrightarrow
        f^\ell x\in Y \text{ for some } 1 \le \ell \le k - 1
        .
      \end{equation}

      \item[A4\label{itm:tail}]\textbf{Regularly varying tails:}
    Let $\tau : Y \to \mathbb{N}$ and $\tau^{(1)}, \ldots, \tau^{ (d) } : Y \to \mathbb{N}$ be given by
      \begin{equation}\label{eq:def-tau}
        \tau (x) \coloneqq \min \{ n \ge 1 :  f^n x \in Y \},
        \text{ and }
        \tau^{(i)} (x) = \Card \{ 1 \le n \le \tau (x) \colon f^n x \in X_i \},\;
        \text{for } i = 1, \ldots, d
        ,
      \end{equation}
      and note that $\tau$ is the first return time to $Y$.
      There are functions $\ell_1,\ldots, \ell_d \colon \R \to \R$ which are slowly varying at $\infty$, so that 
      \begin{equation*}
               \mu ( \tau^{(i)} > n ) \sim \ell_i(n) n^{-1}
        \quand
        \sum_{ n = 1 }^{ \infty } \ell_i (n) n^{ -1 } = \infty.
      \end{equation*}
      Moreover, the limit
      \begin{equation}\label{eq:not-crazy-ells}
        \rho
        _i \coloneqq \lim_{ n \to \infty } \frac{ \ell_i (n) }{ \ell_1 (n) + \cdots + \ell_d (n) },
      \end{equation}
      exists for each $i = 1, \ldots, d$.

      \item[A5\label{itm:mixing}] \textbf{Return times are continued fraction mixing:} Let $F \coloneqq f^{ \tau } \colon Y \to Y$ denote the first return map to $Y$. We assume that $\mu_Y \coloneqq \mu|_Y$ is invariant and ergodic for $F$. Moreover, we assume for each $i = 1,\ldots, d$ the sequence $Z_n = \tau^{(i)} \circ F^{n-1}$ is \emph{continued fraction mixing}
      \begin{equation}\label{eq:def-cf-mixing}
        \sup_n \psi (n) < \infty, \text{ and } \psi(n) \to  0,
        \text{ where }
        \psi ( n )
        \coloneqq 
        \sup_{ k \ge 1 }
        \sup_{ \substack{ A \in \sigma ( Z_1, \ldots, Z_{k}), \\ B \in \sigma ( Z_{k+n}, Z_{k+n+1} \ldots ) \\ \mu_Y ( A )\mu_Y( B ) > 0.  }}
        \left|
          \frac{ \mu_Y ( A \cap B ) }{ \mu_Y (A) \mu_Y (B) }  - 1
        \right|,
      \end{equation}
      and the continued fraction mixing coefficients satisfy $\sum_{ n = 1 }^{ \infty } \psi ( n ) / n < \infty$.
    \end{description}

\begin{remark}\label{rem:dyn-sep}
  Notice that $\tau = \tau^{(1)} + \cdots + \tau^{(d)} + 1$
  and the dynamical separation property of $Y$ in \ref{itm:dyn-sep} implies that
  \[
    \tau^{(i)}  > 0 \Longrightarrow \tau^{(j) } = 0, \text{ for all } i,j = 1, \ldots, d, \text{ with } i \neq j.
  \]
  Moreover, \ref{itm:tail} implies that
  \[
    \mu_Y ( \tau > n )  \sim \sum_{ i = 1 }^{ d }  \mu_Y ( \tau^{ (i) } > n )\sim \ell(n) n^{-1}
    \quad\text{where}\quad
    \ell \coloneqq \ell_1 + \cdots + \ell_d.
  \]
\end{remark}

We define the quantities $L_i, \NN^{(i)}$ in terms of the $\ell_i$ as in \eqref{eq:def-Li-Ni}. We set
\[
  \vrho = ( \rho_1, \ldots, \rho_d),
\]
with $\rho_i$, $i=1,\ldots,d$ as in \eqref{eq:not-crazy-ells} and notice that $\vrho \in \cS$. As in the introduction we define $\RR_i$ to be the line which joins $\vrho$ with the $i^{th}$ vertex $\vv_i = ( 0, \ldots,0 , 1, 0, \ldots, 0 )$ of the simplex $\cS$. Also, for $N \ge 1$ and  $\vz \in \{ 1, \ldots, d \}^N$ let $\RR_{ \vz}$ denote the convex hull of the lines $\RR_{ \zeta_1 }, \ldots, \RR_{ \zeta_N }$.
Then our first main theorem in the general setting reads as follows:

\begin{manualtheorem}{A}\label{thm:2-trimming}
  Suppose that $f$ satisfies~\ref{itm:f-cons}--\ref{itm:mixing}. If  $\NN^{(1)} = \ldots = \NN^{(d)} = N$ for some $N \geq 1$ and $\vrho \in \cS$, then,
  \[
    \cL \aseq \pi \left(\bigcup_{ \vi \in \{ 1, \ldots, d \}^N } \RR_{\vi}\right). 
  \]
\end{manualtheorem}

If the slowly varying function is asymptotic to a constant, then the theorem simplifies analogously as in Theorem~\ref{thm:theorem-constant}.

\begin{corollary}\label{cor:constant}
  Suppose that $f$ satisfies~\ref{itm:f-cons}--\ref{itm:mixing} and $\ell_i \sim c_i$, for some constants  $c_1, \ldots, c_d > 0$. Then,
  \[
    \cL \aseq \pi \left(\bigcup_{ i = 1 }^d \RR_{i}\right).
  \]
\end{corollary}

Similarly as before, we can state an analog to Theorem~\ref{thm:phys-meas-case} in the general setting:

\begin{manualtheorem}{B}\label{thm: phys meas}
 Suppose that $f$ satisfies~\ref{itm:f-cons}--\ref{itm:mixing} and that $\ell_i(n)=c_i \log n$ for $i=1,\ldots, j$ and $\ell_i(n)= c_i$ for $i=j+1,\ldots, d$ for some $c_1, \ldots, c_d > 0$.
 Then,
  \begin{equation*}
  \cL \aseq \pi \left(\bigcup_{ i=1}^j \RR_i\right).
\end{equation*}
\end{manualtheorem}

\begin{remark}
 We note that if $j=1$, i.e.\ there is only one fixed point where the slowly varying function is large, then we obtain a physical measure.
 
 Moreover, the theorems do not necessarily describe \emph{all} of the different possible limit sets $\cL$ that can appear in the setting of \ref{itm:f-cons} to \ref{itm:mixing}. Instead, the last theorem mainly aims to illustrate that one may construct cases where all $\tau^{(i)}$ are non-integrable, but only some $i$ actually play a role for the limiting behaviour of the empirical measure.   
\end{remark}

\section{Continued fraction mixing sequences with 1-regularly varying tails and trimmed strong laws of large numbers}
\label{sec:cf mixing}

Throughout this section we assume that $Z_1,Z_2,\ldots $ is a sequence of identically distributed non-negative random variables on a probability space $(\Omega, \cB, \P)$ so that $\E ( Z_1 ) = \infty$ and assume moreover that the sequence $Z_n$ is \emph{continued fraction mixing} in the sense of \eqref{eq:def-cf-mixing}.
In this section we will collect some results on the almost sure behaviour of
\[
  S_n
  \coloneqq
  \sum_{ k = 1 }^{ n } Z_k,
  \quand
  M_n \coloneqq
  \max_{ 1 \le k \le n } Z_k.
\]

As $\mathbb{E} ( Z_1 ) = \infty$ we cannot expect a strong law of large numbers to hold, indeed for any sequence $d_n \to \infty$ one has that either
\begin{equation}
  \liminf_{ n \to \infty } \frac{S_n}{ d_n } \aseq 0
  \quad\text{or}\quad
  \limsup_{ n \to \infty } \frac{S_n}{ d_n } \aseq \infty,\label{eq: Aaro}
\end{equation}
(see~\cite{ChoRob1961} for a proof in the iid setting and see~\cite{Aar77} in the ergodic setting).
However, a weak law of large numbers may still hold. Let
\begin{equation}\label{eq:def-L-and-b}
  L (x) \coloneqq \int_{ 0 }^x \P ( Z_1 > t ) \, dt
  \quand
  b_n \coloneqq \inf \left\{ t > 0 \colon \frac{ t }{ L(t) } > n \right\}.  
\end{equation}
We then have the following theorem which is due to Feller~\cite{feller} in the iid case and can be found in \cite{Sze2001} for the continued fraction mixing case.\footnote{Note that~\cite[Theorem 3]{Sze2001} requires only $\phi$-mixing which is slightly weaker than our present assumptions.}

\begin{theorem}[{\cite[Theorem 3]{Sze2001}}]\label{thm:weak-law}
  $\displaystyle \frac{ S_n }{ b_n } \to^\P 1$ if and only if $L$ is slowly varying.
\end{theorem}
Moreover from~\cite{Sze2001}, one has the following.\footnote{Note that also in~\cite[Theorem 4]{Sze2001} only $\phi$-mixing is required.}
\begin{theorem}[{\cite[Theorem 4]{Sze2001}}]\label{thm:inf-as}
  If $L$ is slowly varying, then
  $\displaystyle \limsup_{n\to\infty} \frac{ S_n }{ b_n } \aseq \infty$.
\end{theorem}
It turns out that one can sometimes still obtain a strong law for $S_n$ if one trims a certain finite number of maximal elements. The number of maximal elements one has to remove is determined by the quantity
\begin{equation*}
  \NN \coloneqq \min \left\{ k \ge 1 \colon \sum_{ n = 1 }^{ \infty }
      \left( n \frac{ L' (n) }{ L (n) }  \right)^{ k + 1 } \frac{ 1 }{ n }
      < \infty
    \right\},
\end{equation*}
where here and throughout we will make the convention $\min\emptyset = \infty$.
We write $M_{ n,k }$ for the $k^{ th }$ largest term in the sequence $Z_{ 1 }, \ldots, Z_{ n }$.
Then, the following result of Aaronson-Nakada \cite{AarNak2003} provides a \emph{trimmed} strong law of large numbers, provided that the mixing coefficients $\psi(n)$ decay sufficiently quickly.

\begin{theorem}[{\cite[Theorem 1.1]{AarNak2003}}]\label{thm:aar-nak}
  Suppose that $\sum_{ n = 1 }^{ \infty } \frac{ \psi (n) }{ n } < \infty$. Then the following two statements hold:
  \begin{enumerate}
    \item $\displaystyle \limsup_{ n \to\infty} \sum_{ k = 1 }^{ n } 1_{ \{ Z_k > t b_n \} } \aseq \mathfrak{N}$, for every $t > 0$.

    \item If $\mathfrak{N} < \infty$, then
    $S_n \aseq M_{ n,1 } + \cdots + M_{ n, \mathfrak{N} } + b_n + o(b_n)$.
  \end{enumerate}
\end{theorem}

\begin{remark}
  If $\mathfrak{N} < \infty$, then $L$ is slowly varying.
  This can be seen as follows: Define 
  \begin{align}
   \epsilon(t)=t(\log (L))'(t)=\frac{t\mathbb{P}(Z>t)}{L(t)}.\label{eq: def epsilon}
  \end{align}

  Note that $\log (L(wt))-\log ( L(t))=\int_1^w(\epsilon(tu)/u)du$, for any $w>1$.  
  Hence, $L$ is slowly varying if and only if $\lim_{t\to\infty}\int_1^w(\epsilon(tu)/u)du=0$ for all $w>1$. 
  Since $\epsilon(t)/t$ is non-increasing the last statement is equivalent to $\lim_{t\to\infty} \epsilon(t)=0$.

  On the other hand, if for some $k\in\mathbb{N}_0$, we have $\sum_{ n = 1 }^{ \infty }
  \left( n  L' (n) / L (n)   \right)^{ k + 1 } / n \sim \sum_{n=1}^{\infty} (\epsilon(n)^{ k  +  1 } /n)<\infty$, then $\lim_{t\to\infty}\int_1^w(\epsilon(tu)^{ k + 1 } /u)du=0$ and  $\lim_{t\to\infty} \epsilon(t)=0$.
\end{remark}

In what follows we will be interested in sequences whose tail distributions are $1$-regularly varying, i.e.\ 
\begin{equation}\label{eq: PX>t}
   \P ( Z > t ) = \ell (t) t^{-1},
\end{equation}
where $\ell$ is slowly varying. In this case the function $L$ is slowly varying (see e.g.\ \cite[Proposition~1.5.9a]{BinGolTeu1987}) and the expression for $\NN$ becomes
\begin{equation}\label{eq:def-N}
    \NN \coloneqq \min \left\{ k \ge 1 \colon \sum_{ n = 1 }^{ \infty }
      \left( \frac{n \ell (n) }{ L (n) }  \right)^{ k + 1 } \frac{ 1 }{ n }
      < \infty
    \right\}.
\end{equation}
In following example we will give some consequences of Theorems~\ref{thm:weak-law},~\ref{thm:inf-as} and~\ref{thm:aar-nak} in the case that the slowly varying function $\ell$ is asymptotically constant.

\begin{example}[Continued fraction mixing sequences with $\ell \sim c$.]\label{ex: ellsimc}
  Assume that the $Z_n$ have distribution of the form~\eqref{eq: PX>t} with $\ell(t)\sim c$ for some $c > 0$ and that their continued fraction mixing coefficients satisfy $\sum_{ n = 1 }^{ \infty } \psi ( n ) / n < \infty$.
  Inserting this assumption into~\eqref{eq:def-L-and-b}, we find
  \[
    L(t)
    \sim \int_{1}^{t} \frac{c}{s} \, ds
    \sim c \log(t)
  \]
  and hence,
  \[
    b_n \sim c n \log (n),
    \quand
    \NN = 1.
  \]

  So, in this setting Theorem~\ref{thm:aar-nak} tells us that,
  \begin{equation}\label{eq:one-trimming}
    S_n \aseq M_n + c n \log n  + o ( c n \log n )
    \quand
    \limsup_{n\to\infty} \sum_{ k = 1 }^{ n } 1_{ \{ Z_k > t c n \log n \} } \aseq 1
    \quad \text{for all} \quad t > 0.
  \end{equation}
  Combining the first part of~\eqref{eq:one-trimming} with Theorem~\ref{thm:weak-law} gives that
  \begin{equation}\label{eq:Mn-to-prob}
    \frac{ M_n }{ b_n } \aseq \frac{ S_n -  b_n }{ b_n}  \to^\P 0,
  \end{equation}
  and~\eqref{eq:Mn-to-prob} together with Theorem~\ref{thm:inf-as} tells us that
  \begin{equation*}
    \liminf_{n\to\infty} \frac{ M_n }{ b_n } \aseq 0
    \quand
    \limsup_{n\to\infty} \frac{ M_n }{ b_n } \aseq \infty.
  \end{equation*}
\end{example}

\section{Proofs of Theorems \ref{thm:2-trimming} and \ref{thm: phys meas}}\label{sec:proofs}

Throughout this section we will assume that $f \colon X \to X$ is a non-singular transformation of a standard probability space $(X,\cB, m)$ which satisfies \ref{itm:f-cons}--\ref{itm:mixing}.

\subsection{Encoding $\cL$}\label{subsec:coding}

Let $\vS_n =  ( S_n^{(1)}, \ldots, S_n^{(d)})$ where $S_n^{(i)} \coloneqq \sum_{ k = 0 }^{ n - 1 } 1_{ X_i } \circ f^k$.
From \ref{itm:f-cons} and \ref{itm:fxd-pnts} we know by the ergodic theorem for infinite measures (see \cite[Exercise 2.2.1]{Aar1997}) that for any $K \subset X$ which is compact and doesn't contain any of the neutral fixed points one has that
\[
   \lim_{ n \to \infty }   \frac{ 1 }{ n } \sum_{ k = 0 }^{  n - 1 }  1_{ K } \circ f^k \aseq 0.
\]
Thus, $e_n (x) (K) \to 0$ almost surely. From this one can deduce that, almost surely, the limit points of $e_n$ must be contained in $\pi ( \cS )$ (see \cite[Lemma 3.2]{CoaMelTal2026} for details). Once this is established it is not hard to deduce that the limit points of the occupation frequencies $\frac{ 1 }{ n } \vS_n$ necessarily determine the corresponding sequence of empirical measures (see \cite[Corollary 3.3]{CoaMelTal2026}). Hence we have:

\begin{proposition}[{\cite[Lemma 3.2 and Corollary 3.3]{CoaMelTal2026}}]\label{prop:coding}
  $\cL \aseq \pi \left\{\text{limit points of } \frac{ S_n }{ n }\right\} \subset \pi\cS$.
\end{proposition}

\subsection{$\frac{1}{n} \vS_{ n }$ along return times}\label{subsec: sumalongreturns}

Proposition~\ref{prop:coding} reduces understanding the limit points of $e_n$ to understanding the limit points of the occupation times $\frac{1}{m}\vS_m$. It turns out that it is convenient to first understand the behaviour of the $\frac{1}{m}\vS_m$ along the subsequence of times which correspond to returns to $Y$. To that end we let $\vtau \coloneqq ( \tau^{(1)}, \ldots, \tau^{(d)} )$, and $\vtau_n = \sum_{ k = 0 }^{ n - 1 } \vtau \circ F^k$ and note the correspondence
\begin{equation}\label{eq:subsequence}
  \frac{ \vtau_n }{ \tau_n } = \frac{\vS_{\tau_n}}{\tau_n}
  .
\end{equation}

The following lemma will allow us to deduce the limit points of the $\frac{ 1 }{ m }\vS_m$ from the $\frac{ \vtau_n }{ \tau_n }$.

\begin{lemma}\label{lem:dist-Sn-taun}
  For any $\eps > 0$ and for almost every $x \in Y$ there exists an $N \ge 1 $ such that
  \begin{equation}\label{eq:S_n-dist}
    \max_{ \tau_n (x) \le j  \le \tau_{n + 1 } (x) }
    \mathrm{dist} \left(
      \frac{ 1 }{ j }\vS_{ j }(x),
      \left\{
        t\frac{\vtau_n}{\tau_n}(x)
        + ( 1 - t )\frac{\vtau_{n+1}}{ \tau_{n+1}}(x)
        \colon t \in [0,1]
      \right\}
    \right) < \eps,
  \end{equation}
  for all $n \ge N$, where ``$\mathrm{dist}$'' denotes the Hausdorff distance.
\end{lemma}

\begin{proof}
  Let $x \in Y$. From \eqref{eq:subsequence} it suffices to consider $n$ and $j$ such that  $\tau_n ( x ) < j < \tau_{ n + 1 } (x)$. Writing $j = \tau_n ( x ) + k$ for some $1 \le k \le \tau \circ F^n(x) - 1$, we know from \ref{itm:dyn-sep} that $\tau \circ F^n (x) = \tau^{(i)} \circ F^n(x)+1$ for some $i=1,\ldots,d$ and so
  \[
     \frac{ \vS_{ \tau_n(x) + k }(x) }{ \tau_n(x) + k }
     = \left(
      \frac{ \tau_n^{(1)}(x) }{ \tau_n(x) + k } , \ldots, \frac{ \tau_n^{(i)}(x)  +  k }{ \tau_n(x) + k } , \ldots, \frac{ \tau_n^{(d)}(x) }{ \tau_n(x) + k } 
    \right).
  \]
  This means that for almost every $x \in Y$, eventually the set $\{ \vS_{ \tau_n(x) + k }(x) /  (\tau_n(x) + k ) \colon 1 \le k  \le \tau\circ F^n(x) \}$ must be contained in an arbitrarily small neighbourhood of the straight line\footnote{For fixed $x$ and fixed $i$ consider the curve $g : [0, \tau\circ F^n(x) ] \to \R^d$, given by $g(t) = \left( \frac{ \tau_n^{(1)} (x) }{ \tau_n(x) + t } , \ldots, \frac{ \tau_n^{(i)}(x)  +  t }{ \tau_n(x) + t } , \ldots, \frac{ \tau_n^{(d)}(x) }{ \tau_n(x) + t }  \right) $. Then, $g(0) = ( \vtau_n / \tau_n )(x)$ and $g ( \tau\circ F^n(x) ) = ( \vtau_{ n + 1} / \tau_{ n + 1} ) (x)$ and one readily verifies that $g'$ and $g''$ are linearly independent and so $g$ is indeed the straight line joining the two endpoints.} in $\R^{d}$ which joins the points $(\vtau_n / \tau_n)(x)$ and $(\vtau_{n+1}/ \tau_{n+1})(x)$, concluding the proof.
\end{proof}

\subsection{Consequences of the trimmed strong law}\label{subsec:conseq trimmed lln}
  From \ref{itm:tail} and \ref{itm:mixing} we know that the sequences $(\tau \circ F^n), (\tau^{ (1) } \circ F^n), \ldots, (\tau^{ (d) } \circ F^n)$ all satisfy the hypothesis of Theorem~\ref{thm:aar-nak}. Following~\eqref{eq:def-L-and-b} and~\eqref{eq:def-N}, for each $i = 1, \ldots, d$ we let
  \begin{gather*}
    L_i (x) \coloneqq \int_{ 0 }^x \mu_Y ( \tau^{ (i) } > t ) \, dt,
    \quad
    b^{ (i) }_n \coloneqq \inf \left\{ t > 0 \colon \frac{ t }{ L_i(t) } > n \right\},
    \\
    \NN^{(i)} \coloneqq \min \left\{ k \ge 1 \colon \sum_{ n = 1 }^{ \infty }
      \left(\frac{ n \ell_i (n) }{ L_i (n) }\right)^{ k + 1 } \frac{ 1 }{ n }
      < \infty
    \right\},
  \end{gather*}
  and let $M_{ n,k }^{ (i) }$ denote the $k^{ th }$ largest element of $\tau^{ (i) }, \ldots, \tau^{ (i) } \circ F^{ n-1 }$. Theorem~\ref{thm:aar-nak} tells us that if $\NN_i < \infty$ for a given $i = 1, \ldots, d$ and $t>0$, then
  \begin{equation}\label{eq:taui-trimmed-law-1}
        \tau^{ (i) }_n \coloneqq
    \sum_{ k = 0 }^{ n - 1 } \tau^{ (i) } \circ F^{ k } \aseq
    M^{ (i) }_{ n, 1 } + \cdots + M_{ n, \mathfrak{N}^{(i)} }^{(i)} + b_n^{(i)} + o( b_n^{(i)} )
    \quand 
    \limsup_{ n \to \infty }
    \sum_{ k = 0 }^{ n - 1 } 1_{ \{ \tau^{ (i) } \circ F^{ k }  > t b_n^{(i)} \}} \aseq \mathfrak{N}^{(i)}.
  \end{equation} 
  Similarly, let
  \[
    L (x) \coloneqq \int_{ 0 }^x \mu_Y ( \tau > t ) \, dt,
    \;
    b_n \coloneqq \inf \left\{ t > 0 \colon \frac{ t }{ L(t) } > n \right\},
    \;
    \NN \coloneqq \min \left\{ k \ge 1 \colon \sum_{ n = 1 }^{ \infty }
      \left(\frac{n \ell(n) }{ L (n) }\right)^{ k + 1 } \frac{ 1 }{ n }
      < \infty
    \right\},
  \]
  and let $M_{ n,k }$ denote the $k^{ th }$ largest element of $\tau, \ldots, \tau \circ F^{ n-1 }$. Again, Theorem~\ref{thm:aar-nak} tells us that if $\mathfrak{N} < \infty$ and $t>0$ then
  \begin{equation}\label{eq:tau-trimmed-law-1}
    \tau_n \coloneqq
    \sum_{ k = 0 }^{ n - 1 } \tau \circ F^{ k } \aseq
    M_{ n, 1 } + \cdots + M_{ n, \mathfrak{N} } + b_n + o( b_n ),
    \quand
    \limsup_{ n \to \infty }
    \sum_{ k = 0 }^{ n - 1 } 1_{ \{ \tau \circ F^{ k }  > t b_n \}} \aseq \NN.
  \end{equation}
  
The maxima of $\tau$ are necessarily maxima of the $\tau^{ (i) }$ in the sense that for every $k \le n$, 
\begin{equation}\label{eq: max tautaui}
 M_{ n,k } = \tau^{ (i) } \circ F^{ m } + 1 = M_{ n, j }^{ (i) } + 1\text{ for some }i = 1, \ldots, d,\text{ some }m \le n - 1\text{ and some }j \le k.
\end{equation}

\begin{remark}\label{rem:correspondence}
  From \ref{itm:dyn-sep} we know that $\rho_1 + \cdots + \rho_d = 1$. From Remark~\ref{rem:dyn-sep}, we know that
  \[
    L \sim L_1 + \cdots + L_d,
  \]
  and hence
  \[
    \NN = \max \{ \NN_1, \ldots, \NN_d \}, 
  \]
  and
  $b_n^{ (i) } \le b_n$ for every $i = 1, \ldots, d$.
  Moreover~\eqref{eq:not-crazy-ells} ensures that $L_i / L \sim \rho_i$ for each $i = 1, \ldots, d$. Hence $b_n^{(i)} = \rho_i b_n$ whenever $\rho_i > 0$. On the other hand, if $\rho_i = 0$ for some $i = 1, \ldots, d$, then
  \[
    \frac{ b_n^{ (i) } }{ b_n } =
    \frac{ b_n^{ (i) } }{ L_i ( b_n^{ (i) } ) } \frac{ L_i ( b_n^{ (i) } ) }{ L( b_n^{ (i) } ) } \frac{ L ( b_n^{ (i) }) }{ L ( b_n ) } \frac{ L( b_n ) }{ b_n } \ll\frac{ L_i ( b_n^{ (i) } ) }{ L( b_n^{ (i) } ) } \to 0,
  \]
  where in the above we have used the definition of $L$ and $L_i$ to see that $\frac{ b_n^{ (i) } }{ L_i ( b_n^{ (i) } ) } \ll n$ and $\frac{ L( b_n ) }{ b_n } \ll 1 / n$ and we have used Potter's bound (see e.g.\ \cite[Thm.~1.5.6]{BinGolTeu1987}) together with the fact that $b_n^{ (i) } \ll b_n$ to conclude that $\frac{ L ( b_n^{ (i) }) }{ L ( b_n ) } \ll 1$. 
\end{remark}

A significant advantage to studying the sequence $\vtau_n / \tau_n$ over the full sequence $\vS_m / m$ is that a trimmed strong law as in \eqref{eq:taui-trimmed-law-1} holds for the former. We now give some consequences of the trimmed strong law on the possible limit points of $\vtau_n / \tau_n$.
For  $0 \le N \le d$, $\vz  \in \{ 1, \ldots, d \}^{ N }$ and $\vg = (\gamma_1, \ldots, \gamma_{ N } ) \in \R_{ \ge 0 }^{ N }$ we set
\begin{equation}\label{eq:def-R}
  \RR_{ \vz } ( \vg ) \coloneqq
  \frac{ \vrho + \sum_{ i = 1 }^{ N }  \gamma_i \vv_{ \zeta_{ i } } }{ 1 + \gamma_1 + \cdots + \gamma_{ N } }
  ,
\end{equation}
where $N  = 0$ we view both $\vz$ and $\vg$ as the empty tuple $\emptyset$ and set
\[
  \RR_{ \emptyset } ( \emptyset ) = \RR_{ \emptyset } = \vrho.
\]

\begin{lemma}\label{lem:propeties-of-R}
  Let $N \ge 0$, $k \ge 1$, $\vz \in \{ 1, \ldots, d \}^{ N } $, $\vz' = ( \zeta_1, \ldots, \zeta_N, \zeta_{ N + 1 }, \ldots, \zeta_{ N + k } ) \in \{ 1, \ldots, d \}^{ N + k }$, and let $\gamma_1, \ldots, \gamma_{ N + k } \ge 0$. Then, one has
  \begin{enumerate}
    \item\label{itm:prop1} $\bigcup_{ \vg \in \R_{ \ge 0 }^{  N } } \RR_{ \vz } ( \vg ) =  \RR_{ \vz }$.
    
    \item\label{itm:perm} For any $\sigma \colon \R^{ N } \to \R^{ N }$ that is a permutation of coordinates, $\RR_{ \sigma( \vz ) } ( \sigma(\vg) ) = \RR_{ \vz } ( \vg )$.

    \item\label{itm:subset} $\RR_{ \vz' } ( \gamma_1, \ldots, \gamma_{ N }, 0, \ldots, 0 ) = \RR_{ \vz } ( \gamma_1, \ldots, \gamma_N )$, and so  $\RR_{ \vz } \subset \RR_{ \vz' }$.

    \item\label{itm:small} There exists a $C > 0$ so that if $\eps > 0$ and if $\gamma_{ { N+1} }, \ldots, \gamma_{ N + k } < \eps$, then $| \RR_{ \vz' } ( \gamma_1, \ldots, \gamma_{ N + k } ) - \RR_{ \vz } ( \gamma_1, \ldots, \gamma_N ) | < C \eps$.
  \end{enumerate}
\end{lemma}

\begin{proof}
  We first consider Item (\ref{itm:prop1}) in the special case that $N = 1$. Let $\zeta \in \{ 1, \ldots, d \}$. By definition~\eqref{eq:def-R}, $\RR_\zeta (0) = \vrho$ and $\lim_{ \gamma \to \infty } \RR_\zeta (\gamma) = \vv_\zeta$. Moreover, as the derivatives $\frac{ d }{ d \gamma } \RR_\zeta ( \gamma )$ and $\frac{ d^2 }{ d\gamma^2 }\RR_\zeta (\gamma)$ are linearly dependent for each $\gamma$ we know that the curve $\gamma \mapsto \RR_\zeta (\gamma)$ is in fact the straight line that connects $\vrho$ and $\vv_\zeta$. Thus, $\bigcup_{ \gamma \ge 0 } \RR_\zeta  (\gamma) = \RR_\zeta$.

  Now suppose that $N \ge 1$, let $\vz \in \{1, \ldots, d\}^{ N }$, and let $\vg \in \R^{ N }_{ \ge 0 }$. From what we have just shown, for every $i = 1, \ldots, N$, there is a $t_i \in [0,1]$ such that
  \[
    \RR_{ \zeta_i } ( \gamma_i ) = t_i \vrho + (1 - t_i) \vv_{ \zeta_i }.
  \]
  So, 
  \begin{align*}
    \RR_{ \vz } ( \vg ) &=
                        \sum_{ i = 1 }^{ N } \left(\RR_{ \zeta_i } ( \gamma_i ) \frac{ 1 + \gamma_i }{ 1 + \gamma_1 + \cdots + \gamma_N }  - \frac{ N - 1 }{ 1 + \gamma_1 + \cdots + \gamma_N } \vrho \right)\\
                      &=
                         \vrho \sum_{ i = 1 }^{ N }  \frac{  t_i ( 1 + \gamma_i )  - N + 1}{ 1 + \gamma_1 + \cdots + \gamma_N } + \sum_{ j \in \vz } \vv_j \sum_{ i \colon \zeta_i = j } \frac{ (1 - t_i )(1 + \gamma_i ) }{ 1 + \gamma_1 + \cdots + \gamma_N },
  \end{align*}
  from which one readily verifies that $\RR_{ \vz }( \vg )$ is a convex combination of the vectors $\vrho$ and $\{\vv_j \colon j \in \vz\}$, yielding that $\RR_{ \vz } (\vg) \in \RR_{ \vz }$. Moreover, every convex combination of $\vrho$ and $\{\vv_j \colon j \in \vz\}$ can be written in this way. This concludes Item (\ref{itm:prop1}).

  Items (\ref{itm:perm})--(\ref{itm:small}) then follow from Item (\ref{itm:prop1}) and \eqref{eq:def-R}.
\end{proof}

\begin{lemma}\label{lem:taun-over-taun}
  Suppose that $\NN < \infty$. Let $N_n \coloneqq \sum_{ k = 0 }^{ n - 1 } 1_{ \{ \tau \circ F^k > \eps b_n \} }$,
  and let $\vz_n = ( \zeta_{ n, 1 }, \ldots, \zeta_{ n, N_n } ) \in \{1, \ldots, d\}^{ N_n } $ be such that $M_{ n, k } = \tau^{ (\zeta_{n, k }) } \circ F^m$ for some $0 \le m \le n - 1$ for each $k = 1, \ldots, N_n$. Then, there exists a $C > 0$, such that for any $\eps > 0$,
  \begin{equation}\label{eq:taun-tau}
    \left\| \frac{ \vtau_n }{ \tau_n } - \RR_{ \vz_n } \left( \frac{ M_{ n,1 } }{ b_n }, \ldots, \frac{ M_{ n, N_n } }{ b_n } \right) \right\| < C \eps,
  \end{equation}
  almost surely eventually, where here and in the following $\|\cdot\|$ denotes the Euclidean norm. 
\end{lemma}

\begin{proof}
  For $k = 1, \ldots, \NN$ let $\gamma_{ n,k } = \frac{ M_{ n,k } }{ b_n }$.
  We know from~\eqref{eq:tau-trimmed-law-1} that
  \[
    \left| \frac{ \tau_n }{ b_n } - ( 1 + \gamma_{ n,1 } + \cdots + \gamma_{ n, \NN } )
    \right| < \eps,
  \]
  almost surely eventually. We also know from~\eqref{eq:tau-trimmed-law-1}  that almost surely eventually, $N_n \le \NN$. So, as only the first $N_n$ of the $\gamma_{ n, k }$ are possibly larger than $\eps$,
  \[
    \left| \frac{ \tau_n }{ b_n } - ( 1 + \gamma_{ n,1 } + \cdots + \gamma_{ n, N_n } )
    \right| <  (\NN+1)\eps,
  \]
  almost surely eventually.
    
  We can apply a similar argument to each of the $\tau^{ (i) }$. Recall \eqref{eq: max tautaui}.
  By the definition of the $\zeta_{ n,k }$, we know that $N_n^{ (i) } \coloneqq  \sum_{ k = 0 }^{ n - 1 } 1_{ \{ \tau^{ (i) } \circ F^k > \eps b_n \} }= \Card \{ i \colon \zeta_{ n,k } = i, k \leq N_n \}$. Moreover, if $\gamma_{ n,j }^{ (i) } \coloneqq \frac{ M_{ n,j }^{ (i) } }{ b_n } \ge \eps$ for some $j \le n$, then we must have that $\gamma_{ n,j }^{ (i) } = \gamma_{ n, k }$ for some $k \le N_n$. Recall also from Remark~\ref{rem:correspondence} that if $\rho_i > 0$, then $b_n^{ (i) } = \rho_i b_n$, otherwise $b_n^{ (i) } = o( b_n)$. So, using~\eqref{eq:taui-trimmed-law-1}, we obtain
  \[
    \left|  
      \frac{ \tau^{ (i) }_n }{ b_n } -
      \left(
        \rho_i + \sum_{ k \colon \zeta_{ n,k } = i } \gamma_{ n, k }
      \right) 
    \right|
    <  (\NN+1)\eps,
  \]
  almost surely eventually. The result then follows.
\end{proof}
 
\begin{lemma}\label{lem:NN-bar}
  If $\NN < \infty$ then there exist $0 \le \overline{\NN}^{(1)}, \ldots, \overline{\NN}^{(d)} < \infty$ such that for any $\eps > 0$ sufficiently small,
  \begin{equation*}
    \overline{\NN}^{ (i) } \aseq \limsup_{ n \to \infty } \sum_{ k = 0 }^{ n - 1 } 1_{ \{ \tau^{ (i) } \circ F^k > \eps b_n \} },
    \quad\text{for all}\quad i = 1, \ldots, d.
  \end{equation*}
\end{lemma}
\begin{proof}
    Let $i = 1, \ldots, d$.
    If $\rho_i > 0$ then it follows immediately that $\overline{ \NN }^{ (i) } = \NN^{(i)}$.
    Let
    \[
      \varphi (x, \eps ) = \limsup_{ n \to \infty }\sum_{ k = 0 }^{ n - 1 } 1_{ \{ \tau^{ (i) } \circ F^k > \eps b_n \} }(x)
      .
    \]
    For fixed $\eps > 0$, $\varphi ( \cdot, \eps)$ is $F$ invariant and hence, by ergodicity, it is almost surely constant and equal to some $\overline{ \NN }^{ (i) }( \eps )$. As $\NN < \infty$ we know that $\NN^{ (i) } < \infty$ and so  $\overline{ \NN }^{ (i) }( \eps ) \le \NN^{ (i) }$. As $\overline{\NN}^{ (i) } ( \eps )$ increases monotonically as $\eps$ decreases, and as it is necessarily integer valued there must exist an $\eps_0 > 0$ such that $\overline{ \NN }^{ (i) }( \eps )$ is constant for all $\eps < \eps_0$.
\end{proof}

\begin{proposition}\label{prop:one-way}
  Suppose that $\NN < \infty$.
  Then, letting $\overline{\NN}^{(1)}, \ldots, \overline{\NN}^{(d)}$ be as in Lemma \ref{lem:NN-bar} and setting
  \[
    I \coloneqq \{ \vz \in \{ 1, \ldots, d\}^{ \NN} \colon \text{ each index } i \in \{ 1, \ldots, d \} \text{ appears at most } \overline {\NN}^{(i)} \text{ times in } \vz \},
  \]
  we have
  \[
    \cL \subset \pi
    \left(\bigcup_{ \vz \in I } \RR_{ \vz }\right)
  \]
  almost surely.
\end{proposition}
\begin{proof}
  We will first prove the result for $\mu_Y$-typical $x \in Y$.
  Let $\eps > 0$. Let $Y' \subset Y$ be the full measure subset on which~\eqref{eq:S_n-dist},~\eqref{eq:taun-tau} and the conclusion of Lemma \ref{lem:NN-bar} hold. Let $x \in Y'$, and let $n_0 \ge 1$ be large enough so that $N_n^{ (i) } \coloneqq \sum_{ k = 0 }^{ n - 1 } 1_{ \{ \tau \circ F^{ k } > \eps b_n \} } (x) \le \overline{\NN}^{ (i) }$ for all $i = 1, \ldots, d$ and so that $\frac{ b_{ n + 1 } }{ b_n } < 2$  for all $n \ge n_0$. Throughout the remainder of the proof we will write $C$ for a positive constant whose value may change from line to line and which is independent of $x,n$ and $\eps$.
  Let $\vz_n$ and $\vg_n$ be as in Lemma \ref{lem:taun-over-taun} so that for every $n \ge n_0$,
  \begin{equation*}
       \left\|
      \frac{ \vtau_n }{ \tau_n }(x) - \RR_{ \vz_n } (\vg_n)
    \right\|
    < C \eps.
  \end{equation*}

  \begin{lemma}\label{sublem:line}
    For $n \ge n_0$, the line $E$ in $\R^d$ joining the points $\frac{ \vtau_n }{ \tau_n }$ and $\frac{ \vtau_{ n + 1 } }{ \tau_{ n + 1 } }$ has Hausdorff distance to $\bigcup_{ \vz \in I } \RR_{ \vz }$ which is no more than $C\eps$.
  \end{lemma}

  \begin{proof}
    Let $N_n = \sum_{ i = 1 }^{ d } N_n^{ (i) }$. Since $N_n \le \NN$, $\RR_{ \vz_n } \subset \bigcup_{ \vz \in I } \RR_{ \vz }$ for every $n$. We divide the proof into different cases.

    \textbf{Case 1:} Suppose that $\frac{ M_{ n, i } }{ b_{ n + 1 } } > \eps$ for all $i = 1, \ldots, N_n$. If $\tau \circ F^{ n } \le \eps b_{ n + 1 }$, then $N_{ n + 1 } = N_n$ and $\vz_{ n + 1 } = \vz_n$
    So, $\dist (\frac{ \vtau_{ n } }{ \tau_{ n } }, R_{ \vz_n } ) < C \eps$, and $\dist (\frac{ \vtau_{ n + 1 } }{ \tau_{ n + 1 } }, \RR_{ \vz_n } ) < C \eps$ and so $\dist( E, \RR_{ \vz_n } ) < C\eps$ by the convexity, of $\RR_{ \vz_n }$.
    On the other hand, if $\tau \circ F^{ n } > \eps b_{ n + 1 }$, then we must have that $N_{ n + 1 } = N_n + 1$ (and hence $N_n \le \overline{ \NN } - 1$) and $\vz_{ n + 1 }$ contains all of the indices in $\vz_n$ together with one additional element. Hence, by Item~\eqref{itm:subset} of Proposition~\ref{lem:propeties-of-R}, $\RR_{ \vz_n } \subset \RR_{ \vz_{ n + 1 } }$. So, again by convexity, $\dist( E, \RR_{ \vz_{ n + 1 }} ) < { C}\eps$. This concludes the proof in the first case.

    \textbf{Case 2:} Suppose that there is a $1 \le k \le N_n$ such that $\frac{ M_{ n,i } }{ b_{ n + 1 } } \le \eps$ for all $i = N_n - k, \ldots, N_n$. Then, it follows by our choice of $n_0$ that $\frac{ M_{ n,i } }{ b_{ n } } = \frac{ M_{ n, i } }{ b_{ n + 1 } } \frac{ b_{ n + 1 } }{ b_n } < 2 \eps$. Hence, by Item~\eqref{itm:small} of Lemma~\ref{lem:propeties-of-R}, we know that $\dist ( \frac{ \vtau_n }{ \tau_n } , \RR_{\zeta_1, \ldots, \zeta_{ N_n - k }}) < C \eps$. Moreover, $\vz_{ n + 1 }$ must contain the indices  $\zeta_1, \ldots, \zeta_{ N_n - k }$, and by Items~\eqref{itm:perm} and~\eqref{itm:subset} of Proposition~\ref{lem:propeties-of-R}, we know that $\RR_{ \zeta_1, \ldots, \zeta_{ N_n - k } } \subset \RR_{ \vz_{ n + 1 } }$. Then, by convexity $\dist( E, \RR_{ \vz_{ n + 1 } } ) < C\eps$. This concludes the second case and the proof of the lemma.
  \end{proof}

  Applying Lemma~\ref{sublem:line} repeatedly and using~\eqref{eq:S_n-dist} we obtain
  \[
    \dist \left(\frac{ S_k(x) }{ k }, \bigcup_{ \vz \in I } \RR_{ \vz } \right) < C \eps
  \]
  for all $k \ge \tau_{ n_0 } (x)$.

  We can now take a sequence $\eps_{ k } \to 0$, and for each $k$ find a full measure set on which $ \dist \left(\frac{ S_k(x) }{ k }, \bigcup_{ \vz \in I } \RR_{ \vz } \right) < C \eps_{ k}$ eventually. Intersecting all of these full measure sets we obtain another full measure set on which
  \[
    \left\{ \text{limit points of } \frac{ S_{ k } }{ k } \right\} \subset
    \bigcup_{ \vz \in I } \RR_{ \vz }.
  \]
  This concludes the result for $\mu_Y$-typical $x \in Y$. Finally, by conservativity, we know that $\mu$-typical points in $X$ will enter $Y$ in finite time which, together with Proposition~\ref{prop:coding}, concludes the proof.
 \end{proof}

 \begin{remark}\label{rem:one-way}
   Suppose that $\NN < \infty$ and that $\rho_i > 0$ for each $i = 1, \ldots, d$ then we know from Remark~\ref{rem:correspondence} that $b_n^{(i)} = \rho_i b_n$ and so Theorem~\ref{thm:aar-nak} guarantees that the hypotheses of Lemma~\ref{prop:one-way} are satisfied with $\overline{\NN}^{ (i) } = \NN^{ {i} } = \NN$ for every $i = 1, \ldots, d$. In this case we find that
   \[
    \cL \subset
    \pi \left(\bigcup_{ \vz \in \{1, \ldots, d \}^{ \NN } } \RR_{ \vz }\right),
  \]
  $\mu$ almost everywhere. 
\end{remark}

The trimmed strong law also gives a partial converse to Proposition~\ref{prop:one-way} in the following sense.

\begin{lemma}\label{lem:other-way}
  $\{ \pi( \vrho ) \} \cup \{ \delta_{\xi_i} \colon \rho_i > 0 \} \subset \cL$ almost everywhere.
\end{lemma}

\begin{proof}
  Suppose that $i \in \{ 1, \ldots, d \}$ is such that $\rho_i > 0$. Then, from Remark~\ref{rem:correspondence} we know that $b_n^{ (i) } = \rho_i b_n$ and so as $\NN^{ (i) } \ge 1$ we know that the event $\{ M_{ n, i  } > t b_n \}$ occurs infinitely often. Hence, from Lemma~\ref{lem:taun-over-taun} and Proposition~\ref{prop:coding} we find that $\delta_{ \xi_i } \in \cL$.

  As $L$ is slowly varying we know that the weak law in Theorem~\ref{thm:weak-law} also holds. Combing this with~\eqref{eq:tau-trimmed-law-1}, we find that
  \[
    \frac{ M_{ n } }{ b_n } \aseq \frac{ \tau_n - b_n }{ b_n } \to 0,
  \]
  in probability. Hence $\liminf_{n\to\infty} M_n / b_n \aseq 0$. So, using again Lemma~\ref{lem:taun-over-taun} and Proposition~\ref{prop:coding} we find that $\pi(\vrho) \in \cL$.
\end{proof}

As an almost immediate corollary to Proposition~\ref{prop:one-way} and Lemma~\ref{lem:other-way} we obtain Corollary~\ref{cor:constant}.

\begin{proof}[Proof of Corollary~\ref{cor:constant}]
  By assumption $\NN = 1$ and $\rho_i > 0$ for each $i = 1, \ldots, d$ and so by Proposition~\ref{prop:one-way} and Lemma~\ref{lem:other-way} we find that
  \[
    \{ \pi( \vrho ) \} \cup \{ \delta_{\xi_i} \colon \rho_i > 0 \}
    \subset \cL \subset \pi \left( \bigcup_{ \vz \in \{1, \ldots, d \}^{ \NN } } \RR_{ \vz } \right).
  \]
  The result follows by the fact that $\cL$ is connected.
\end{proof}

\subsection{Proof of Theorem \ref{thm:2-trimming}}\label{subsec: ProofthmA}
In this section we will assume, in addition to \ref{itm:f-cons}--\ref{itm:mixing}, that $\NN^{ (1) } = \cdots = \NN^{ (d) } = \NN < \infty$.

We first state the main lemma that is needed to prove Theorem~\ref{thm:2-trimming}.
\begin{lemma}\label{lem:2-trimming}
  Let $\zeta_{ n, 1 }, \ldots, \zeta_{ n, \mathfrak{N} } \in \{1, \ldots, d\} $ be such that $M_{ n, k } = \tau^{ (\zeta_{n, k }) } \circ F^m$ for some $0 \le m \le n - 1$ for each $k = 1,\ldots, \mathfrak{N}$. Then, for every
  $\gamma_1 \geq \ldots\geq \gamma_{\mathfrak{N}} > \eps > 0$ and $(\zeta_1, \ldots, \zeta_{\mathfrak{N}}) \in \{ 1, \ldots, d \}^{\mathfrak{N}}$, the event
  \[
    B_n \coloneqq \{ \zeta_{ n, 1 } = \zeta_1 ,\ldots, \zeta_{ n, \mathfrak{N} } = \zeta_{\mathfrak{N}}, \quad M_{ n,1 } / b_n \in [\gamma_1 - \eps, \gamma_1 + \eps], \ldots, M_{ n, \mathfrak{N}}  / b_n \in [\gamma_{\mathfrak{N}} - \eps, \gamma_{\mathfrak{N}} + \eps] \},
  \]
  occurs infinitely often almost surely.
\end{lemma}

Before proving this lemma let us show how it implies Theorem~\ref{thm:2-trimming}.

\begin{proof}[Proof of Theorem~\ref{thm:2-trimming}]
  Lemma~\ref{lem:2-trimming} together with Lemma~\ref{lem:taun-over-taun} give that
  \[
    \pi \left( \bigcup_{ \vz\in \{1,\ldots, d\}^{ \NN }} \RR_{ \vz } \right) \subset \cL.
  \]
  Proposition~\ref{prop:one-way} (and in particular Remark~\ref{rem:one-way}) give the other inclusion and conclude the theorem.
\end{proof}

\begin{remark}\label{rem:indep}
  There exists a $\kappa \ge 1$ such that for any sequence $\varphi_n = \varphi \circ \vtau \circ F^n$ with $\varphi \colon \N^d \to \R$ measurable,
  \[
    \mu_Y ( A \cap B ) \approx \mu_Y (A) \mu_Y (B),
    \quad
    \text{for all}
    \quad
    A \in \sigma ( \varphi_{ 1 } , \ldots, \varphi_{ k } ),
    \quad
    B \in \sigma ( \varphi_{ k + \kappa }, \varphi_{ k + \kappa + 1 }, \ldots ),
    \quad
    k \ge 1.
  \]
\end{remark}

We will also make use of the following version of the second Boreli-Cantelli lemma, originally due to R\'enyi and given in its form below in \cite{Pet2002}.
\begin{lemma}[{\cite[Theorem 2.1]{Pet2002}}]
  \label{lem:renyi}
  If $A_n \subset Y$ is a sequence of positive measure events so that $\mu_Y ( A_n \cap A_k ) \ll \mu_Y ( A_n ) \mu_Y ( A_k )$ for every $n \neq k$, then $\mu_Y( \limsup_{n\to\infty} A_n ) > 0$ if and only if $\sum_{ n = 1 }^{ \infty } \mu_Y ( A_n ) = \infty$.
\end{lemma}
\begin{proof}[Proof of Lemma~\ref{lem:2-trimming}]
  Fix $\gamma_1 \geq \ldots\geq \gamma_{\NN} > \eps > 0$ and $(\zeta_1, \ldots, \zeta_{\NN}) \in \{ 1, \ldots, d \}^{\NN}$, and for each $j = 1, \ldots, \NN$ let 
  \begin{equation*}
    C_{ n,j } = \left\{ (\gamma_j - \eps) b_{ n} \le \tau^{(\zeta_j)} \le ( \gamma_j + \eps ) b_{ n } \right\}
  \end{equation*}
  and set
  \[
    A_n  =\bigcap_{j=1}^{\mathfrak{N}} \bigcup_{k\in J_{n,j}} F^{-k}C_{n,j}
    \quwhere
    J_{n,j}=\left[\frac{n}{3}\left(2+\frac{2j-2}{2\mathfrak{N}-1}\right), \frac{n}{3}\left(
        2 + \frac{2j - 1}{ 2\mathfrak{N}-1 }
      \right)\right] \cap \mathbb{N}.
  \]

In words, the event $A_n$ describes the situation that $\mathfrak{N}$ large (relative to $b_n$) events occur between time $2n/3$ and $n$, each of which in a distinct time interval $J_{n,j}$.
We know from \eqref{eq:tau-trimmed-law-1}, that $\limsup_{ n \to \infty } \sum_{ k = 1 }^{ n - 1 } 1_{ \{\tau \circ F^k > (\gamma_{\NN} - \eps) b_n \} } \aseq \NN$ and so almost surely eventually there can be no more than $\NN$ large events. Hence,
$$\limsup_{n\to\infty} A_n \cap Y^*\subset \limsup_{n\to\infty}B_n \cap Y^*$$
for some $Y^* \subset Y$ with $\mu_Y(Y^*)=1$. 
We will show that $A_n$ occurs infinitely often almost surely.

Note that the time intervlas $J_{n,j}$ are spaced $O(n)$ appart and hence we know from Remark~\ref{rem:indep} that
\begin{align}
 \mu_Y (  A_n )
 \gg \prod_{j=1}^{\mathfrak{N}}\mu_Y\bigg(\bigcup_{k\in J_{n,j}} F^{-k}C_{n,j}\bigg)
   \label{eq: A' estim}
\end{align}
for $n$ sufficiently large. Next, we will analyse the single factors independently. The next estimates use ideas from \cite[p.~177/178]{AarNak2003}. 
First we estimate for any $U>1$
\begin{align*}
  \mu_Y(C_{n,j})
  &\ge \mu_Y\left(\tau^{(\zeta_j)}> (\gamma_j - \eps) b_{ n }\right)- \mu_Y\left(\tau^{(\zeta_j)} > ( \gamma_j + \eps ) b_{ n } \right)\\
  &= \frac{\ell ((\gamma_j-\eps) b_{n})}{(\gamma_j-\eps) b_{n}}- \frac{\ell ((\gamma_j+\eps) b_{n})}{(\gamma_j+\eps) b_{n}}\\
  &=\frac{\ell ((\gamma_j-\eps) b_{n})(\gamma_j+\eps)- \ell ((\gamma_j+\eps) b_{n})(\gamma_j-\eps)}{(\gamma_j-\eps)(\gamma_j+\eps) b_{n}}\\
  &\ge \frac{U^{-1} \ell ( b_{n})(\gamma_j+\eps)- U\ell ( b_{n})(\gamma_j-\eps)}{(\gamma_j-\eps)(\gamma_j+\eps) b_{n}},
\end{align*}
for $n$ sufficiently large. Here, we have used Potter's bound for slowly varying functions, see e.g.\ \cite[Thm.~1.5.6]{BinGolTeu1987} and used that we can choose constants $V>1$, $\rho>0$ such that $V^{-1}\left(\frac{1}{\gamma_j-\eps}\right)^{\rho}> U^{-1}$ and 
$V\left(\gamma_j+\eps\right)^{\rho}< U$.
As $U$ was arbitrarily close to $1$, we can choose it so small that $U^{-1} (\gamma_j+\eps)- U(\gamma_j-\eps)<\eps$
and hence $\mu_Y(C_{n,j}) \ll b_n^{-1} \ell (b_n)$. Thus,
\begin{equation}
 \mu_Y(C_{n,j}) 
 \approx \frac{\ell(b_{n})}{b_{n}}\label{eq: Cn estim}
\end{equation}

Let $\kappa \in \N$ be as in Remark \ref{rem:indep}, let 
\[
  D_{n,k,j}
  =  F^{-k}C_{n,j}
  \cap
  \left[
    \bigcap_{u\in J_{n,j} \setminus [k- \kappa, k+ \kappa]} \{\tau\circ F^{u}\le \eps b_n\}
  \right],
\]
and notice that
\begin{equation}\label{eq:dnkj}
  \sum_{k\in J_{n,j}}1_{D_{k,j}}\le 2 \NN\,
    1_{ \bigcup_{ k \in J_{n,j} }  F^{ -k } C_{n,j} }.
\end{equation}

Next, assume that $\lim_{n\to\infty} n \mu_Y(\tau>b_n)=0$ and
assume that $n$ is large enough such that $n \mu_Y(\tau>b_n)<\frac{1}{2}$. 
Then, using Remark \ref{rem:indep}, we have 
\begin{align}
 \mu_Y(D_{ n,k,j })
 &\gg \mu_Y\left(\bigcap_{u\in J_{n,j} \setminus [k-\kappa, k+ \kappa]}\{\tau\circ F^{u}< \eps b_n\} \right)\,\mu_Y( F^{-k}C_{n,j})\notag\\
  &\gg \big(1- n \mu_Y( \tau> \eps b_n) \big)\,\mu_Y(C_{n,j})\notag\\
  & \gg  \mu_Y(C_{n,j}),\label{eq: muDk}
\end{align}
for $n$ sufficiently large.
Hence, integrating \eqref{eq:dnkj} and using \eqref{eq: muDk} and \eqref{eq: Cn estim}, we obtain
\begin{align*}
  \mu_Y\bigg(\bigcup_{k \in J_{n,j}} F^{-k}C_{n,j}\bigg)
 \gg \sum_{k\in J_{n,j}}\mu_Y(D_{n,k,j})
 \gg n\mu_Y(C_{n,j})
 \gg \frac{n\, \ell(b_n)}{b_n}.
  \end{align*}
  From \eqref{eq: A' estim}, we then obtain
\begin{equation}\label{eq:An}
 \mu_Y(A_n)
 \gg \left( \frac{n\, \ell(b_n)}{b_n} \right)^{\NN}. 
\end{equation}

We will now follow the ideas of \cite[Lemma 3]{Mori76} in order to show that $\mu_Y( \limsup_{ n \to \infty } A_n ) = 1$.
From \eqref{eq:An},
\begin{align*}
 \mu_Y(A_{2^n})
 \gg \left( 2^n\, \mu_Y(\tau> b_{2^n})\right)^{\mathfrak{N}}
 \gg \sum_{v\in [2^{n-1}, 2^n]} v^{\mathfrak{N}-1} \mu_Y(\tau> b_{v})^{\mathfrak{N}}.  
\end{align*}
By \cite[Proof of (i) of Theorem 1.1.]{AarNak2003} it follows that 
$\sum_{v=1}^{\infty} v^{\mathfrak{N}-1} \mu_Y(\tau> b_{v})^{\mathfrak{N}}\approx \sum_{v=1}^{\infty} \frac{ \epsilon (v)^{\mathfrak{N}}}{ v }=\infty$
where $\epsilon$ as in \eqref{eq: def epsilon}.
In particular, we have that $\sum_{n=1}^{\infty} \mu_Y ( A_{ 2^n } ) =\infty$. 
Moreover, the event $A_{ 2^n }$ only considers $\tau\circ F^{\lceil 2^{n+1}/3\rceil-1},\ldots, \tau\circ F^{2^{n  } - 1}$ and hence, by Remark \ref{rem:indep},
$\mu_Y ( A_{ 2^n } \cap A_{ 2^k } ) \ll \mu_Y ( A_{ 2^n } )\mu_Y (  A_{ 2^k })$ whenever $k \neq n$. Hence we can apply Lemma~\ref{lem:renyi} to obtain that $\mu_Y( \limsup_{ n \to \infty } A_n ) > 0$. Finally, as the set $\limsup A_n$ is invariant for $F$ we know by the ergodicity of $\mu_Y$ that  $\mu_Y( \limsup_{ n \to \infty } A_n ) = 1$ as required. This concludes the lemma in the case that $n \mu_Y ( \tau > n ) \to 0$.

Now suppose that $\limsup_{ n \to \infty } n \mu_Y ( \tau > n ) > 0$.
Then there exist a subsequence $(n_u)$ and $\delta>0$ such that $n_u \mu_Y(C_{n_u,j})>\delta$. On this subsequence, we note that 
\begin{align*}
 \mu_Y\bigg(\bigcup_{k\in J_{n_u,j}} F^{-k}C_{n_u,j}\bigg)\geq 
 \mu_Y\bigg(\bigcup_{k\in \tilde J_{n_u,j}} F^{-k}C_{n_u,j}\bigg),
\end{align*}
where we choose
$$\tilde{J}_{n_u,j}=\left[\frac{n_u}{3}\left(2+\frac{2j-2}{2\mathfrak{N}-1}\right), \frac{\tilde{n}_u}{3}\left(2+\frac{2j-2}{2\mathfrak{N}-1}\right)\right]\cap \mathbb{N}=:\{n_{n,j},\ldots, \tilde{n}_{u,j}\}$$ 
with $\tilde{n}_{u,j}=\max\{w\in\mathbb{N}:(w-n_{n,j})\mu_Y(\tau>b_{n_u})< 1/2\}$ and define 
\[
  \widetilde{D}_{n,k,j}=
  F^{-k}C_{n, j}\cap \bigcap_{u\in \tilde J_{n,j} \setminus [k-\kappa, k+\kappa]} \{\tau\circ F^{u}\le \eps b_n\}.
\]
Similarly to the previous case we find that
\begin{equation*}
  \sum_{k\in \tilde J_{n,j}}1_{ \widetilde D_{n,k,j}}\le 2 \NN\,
  1_{ \bigcup_{ k \in \tilde J_{n,j} }  F^{ -k } C_{n,j} }.
\end{equation*}
Hence, an analogous calculation to that in \eqref{eq: muDk} together with \eqref{eq: Cn estim} and our choice of $\tilde n_{ n,j }$ yields
\begin{align*}
 \mu_Y\bigg(\bigcup_{ m \in \tilde J_{ n_u,j }} F^{-m}C_{n_u,j}\bigg)\gg 1.
\end{align*}
Following that we can continue as before and use Renyi's Borel-Cantelli Lemma \ref{lem:renyi} (possibly after choosing another subsequence of $(n_u)$) and the ergodicity of $\mu_Y$ to conclude that $\mu_Y( \limsup A_n ) = 1$. This concludes the proof.
\end{proof}

\subsection{Proof of Theorem~\ref{thm: phys meas}}\label{subsec: ProofthmB}
We now assume that in addition to \ref{itm:f-cons}--\ref{itm:mixing} we also have that there are constants $c_1, \ldots, c_d > 0$ and a $1 \le j \le d - 1$ such that $\ell_i (n) = c_i \log (n)$ for $i = 1, \ldots, j$ and $\ell_i \equiv c_i$ for $i = j + 1, \ldots, d$. 

We first start with the following lemma:
\begin{lemma}\label{lem: main lem phys meas}
  Lemma \ref{lem:NN-bar} holds with $\overline{\NN}^{ (1) } = \cdots = \overline{\NN}^{ (j) } = 1$ and let $\overline{\NN}^{ (j + 1 ) }, \ldots, \overline{\NN}^{ (d) } = 0$.   
\end{lemma}

\begin{proof}
First we look at $i\in \{1,\ldots, j\}$. 
We may calculate
  \[
    L_i (t) 
    = c_i \int_{0}^{t} \frac{\log(s)}{s} ( 1 + o(1) ) \, ds
    \sim  c_i  \frac{(\log(t))^2}{2}.
  \]
  So,
  $b_n^{ (i) } = c_i\frac{1}{2} n (\log (n))^2$.
  On the other hand, by Example \ref{ex: ellsimc} we already have $b_n^{(i)}=c_in\log(n)$ if $i\in \{j+1,\ldots,d\}$. 
  Hence, we can conclude that $\sum_{i=1}^d b_n^{(i)}\sim \frac{j}{2}n(\log(n))^2$. Moreover, for $i\in \{1,\ldots, j\}$ we have
  $t   L_i' (t) / L_i (t)    \sim 2/\log(t)$
  and thus $\mathfrak{N}^{(i)} = \overline{\mathfrak{N}}^{(i)}=1$.

In order to calculate $\overline{\mathfrak{N}}^{(i)}$ for $i\in \{j+1,\ldots, d\}$
note that 
for all $\varepsilon>0$ we have
\begin{align*}
 \bigcup_{k=2^n}^{2^{n+1}} \{M_k^{(i)}>\varepsilon k(\log(k))^2\}
 \subset \{M_{2^{n+1}}^{(i)}>\varepsilon 2^n(\log(2^n))^2 \}.
\end{align*}
Moreover, 
\begin{align*}
 \mu_Y\left(M_{2^{n+1}}^{(i)}>\varepsilon 2^n(\log(2^n))^2 \right)
 &\leq 2^{n+1} \mu_Y\left(\tau^{(i)}>\varepsilon 2^n(\log(2^n))^2\right)
 \ll \frac{2^{n+1}}{2^n(\log(2^n))^2}
 \ll \frac{1}{n^2}. 
\end{align*}
As the last term is summable, we may conclude by the first Borel-Cantelli lemma that $\overline{\mathfrak{N}}^{(i)}=0$ for all $i\in \{j+1,\ldots, d\}$. 
\end{proof}

\begin{proof}[Proof of Theorem~\ref{thm: phys meas}]
Combining Lemma~\ref{lem: main lem phys meas} and Proposition~\ref{prop:one-way} immediately implies $ \cL \subset \pi \left(\bigcup_{ i=1}^j \RR_i\right)$ almost everywhere.

On the other hand, Lemma~\ref{lem: main lem phys meas} implies $\rho_i>0$ for $i\in \{1,\ldots, j\}$ and $\rho_i=0$ fro $i\in \{j+1,\ldots, d\}$. 
Hence, by Lemma~\ref{lem:other-way} we obtain $\{ \pi( \vrho ) \} \cup \bigcup_{i=1}^j\{ \delta_{\xi_i} \} \subset \cL$ and since $\cL$ is connected the statement follows. 
\end{proof}

\section{Proofs of Theorems \ref{thm:theorem-constant}--\ref{thm:parabolic}}

\subsection{Preliminaries on Gibbs-Markov Maps}\label{sec:examples}

We begin with some extra preparation needed to apply Theorems \ref{thm:2-trimming} and \ref{thm: phys meas} to the examples described in Section \ref{sec:intro} and in order to prove Theorems \ref{thm:theorem-constant}, \ref{thm:N}, \ref{thm:phys-meas-case} and \ref{thm:parabolic}. In particular, we will need to verify \ref{itm:f-cons}--\ref{itm:mixing} for both the intermittent interval maps defined in Section \ref{sec:thaler-slow} and the parabolic rational maps considered in Section \ref{sec:parabolic}.
For both classes a key step in the verification of \ref{itm:f-cons}--\ref{itm:mixing} is showing the existence of a \emph{Gibbs-Markov} first return map. Let us briefly recall some standard terminology as well as some well-known properties of maps which admit Gibbs-Markov first return maps. In this section $X$ will always denote a metric space equipped with a Borel probability measure $m$. A standard reference for this material is \cite[Section 4]{Aar1997}.

\begin{definition}
  \label{def:gm}
  We say that a non-singular transformation $f \colon X \to X$ together with a measurable partition $\cP$ of $X$ into open subsets is a \emph{Markov map} if the restriction of $f$ to each partition element is a measurable bijection onto a union of partition elements and if the dynamical refinements of $\cP$ generate the Borel $\sigma$-algebra $\cB$ in the sense that $\bigvee_{ n = 0 }^{ \infty } f^{ -n } \cP = \cB$.

  Let $Y$ be a measurable union of partition elements. We say that the first return map $F \colon Y \to Y$ is \emph{Gibbs-Markov} if $F$ is Markov with respect to some measurable partition $\cQ$ of $Y$ which is a refinement of $\cP$ and if in addition the following properties hold: (1) \emph{finite images}: $F \cQ$ is a finite partition of $Y$; (2) \emph{log-Lipschitz distortion}: there exists a $C > 0$ such that
  \[
    \left|
      \log \frac{ v_a' x }{ v_a' y }
    \right|
    \le C d (x,y),
    \quad
    \text{for all } a \in \cQ, x,y \in F(a),
  \]
  where for each $a \in \cQ$, $v_a \colon Fa \to a$ is given by $v_a = (F|_a)^{-1}$ and
  \[
    v_a' = \frac{ d( m \circ v_a ) }{ dm }.
  \]
\end{definition}

\begin{lemma}[{\cite[Section 4]{Aar1997}}]\label{lem:gm-inv-meas}
  Suppose that $f$ is a topologically mixing Markov map with respect to $\cP$ and that $Y$ is a union of partition elements so that the first return map $F = f^\tau \colon Y \to Y$ is Gibbs-Markov with respect to a partition $\cQ$ of $Y$. Then, $F$ preserves an ergodic probability measure $\mu_Y \ll m|_Y$ whose density $h_Y$ is bounded away from $0$ and $\infty$ and is Lipschitz on the elements of $F\cQ$.
  Moreover, $f$ preserves the $\sigma$-finite measure $\mu$ defined by
  \begin{equation}\label{eq:def-mu}
    \mu ( A ) = \sum_{ n = 0 }^{ \infty } \mu_Y ( f^{-n} A \cap \{ \tau \ge n  + 1 \} ). 
  \end{equation}
\end{lemma}

\begin{remark}\label{rem:inv-measure}
  It follows immediately that the measure $\mu$ defined in \eqref{eq:def-mu} satisfies
  \[
    \mu ( X ) < \infty \text{ if and only if } \int \tau \,d \mu_Y < \infty.
  \]
\end{remark}

\begin{lemma}[{\cite[p.~201]{AarDen2001}}]\label{lem:cf-mix}
  Suppose that $F$ is a topologically mixing Gibbs-Markov map and that $\varphi \colon Y \to \R$ is constant on partition elements. Then, $\varphi \circ F^{ n }$ is exponentially continued fraction mixing.
\end{lemma}

\subsection{Intermittent Interval Maps}\label{subsec: inter-inter-maps}

In this section we prove Theorems \ref{thm:theorem-constant}, \ref{thm:N} and \ref{thm:phys-meas-case}. Let $X = [0,1]$, $m = \Leb$, and let $f \colon X \to X$ be a map satisfying \ref{itm:T1}--\ref{itm:T3} and \ref{itm:exp-slow}. Theorems \ref{thm:theorem-constant}, \ref{thm:N} and \ref{thm:phys-meas-case} follow from Corollary \ref{cor:constant} and Theorems \ref{thm:2-trimming} and \ref{thm: phys meas} provided one can show that $f$ satisfies \ref{itm:f-cons}--\ref{itm:mixing} whenever \eqref{eq:rho_i} and \eqref{eq:ellnonint} hold. 
Hence it suffices to prove the following.

\begin{proposition}\label{prop:thaler}
  If \eqref{eq:rho_i} and \eqref{eq:ellnonint} hold, then the map $f$ satisfies \ref{itm:f-cons}--\ref{itm:mixing}.
\end{proposition}

We follow the argument in \cite[Example 4.15]{CoaMel2025-sub}, making the necessary adaptations to treat the addition of slowly varying functions.

We divide the proof of this result into several steps.
First, from \cite{Zwe1998} and the fact that $f$ is full branched we know that \ref{itm:T1}--\ref{itm:T3} imply that $f$ is conservative and ergodic and preserves a unique (up to scaling) $\sigma$-finite measure $\mu \ll m$, and $\mu ( X ) = \infty$. This gives \ref{itm:f-cons}.

In what follows it will be convenient to label the subintervals in $\cP$, $I_1,I_2,\ldots$ so that $\xi_i \in I_i$ for each $i = 1,\ldots, d$. We write $D = \Card \cP$ so that $D - d \in \N \cup \{\infty\}$ denotes the number of uniformly expanding branches of $f$ and write $f_k = f|_{I_k}$.

Our next step is to construct a suitable decomposition $X = X_1 \cup \cdots \cup X_d \cup Y$.
For $i = 1 ,\ldots, d$ we set $X_i \coloneqq I_i \cap f^{-1}_i I_i$.
The choice of the $Y$ depends on the values of $d$ and $D$. First of all, if $d = D = 2$
there exists a unique periodic orbit $\{ \gamma, f \gamma \}$ of period 2, in this case we set $Y \coloneqq [ \gamma, f \gamma ]$. 
If instead $d = D \ge 3$, set $Y \coloneqq \bigcup_{ i = 1 }^{ d } I_{ i } \setminus f_i^{ -1 } I_{ i }$.  Otherwise, if $1 \le d < D$, and so there is always at least one uniformly expanding branch, we take  $Y = \bigcup_{ i = 1 }^{ d } ( I_i \setminus f_i^{-1}I_i ) \cup \bigcup_{ i = d + 1 }^{ D } I_i$. Notice we always have that $X = X_1 \cup \cdots \cup X_d \cup Y$.
We will now focus on the case that $1 \le d < D$ as, after having defined $Y$, the situation in which $2 \le d = D < \infty$ is analogous.

Define the functions $\tau, \tau^{ (1) }, \ldots, \tau^{ (d) } \colon Y \to \mathbb{N}$ as in~\eqref{eq:def-tau} by setting,
\[
  \tau (x) \coloneqq \min \{ n \ge 1 \colon f^{ n } x \in Y \},
  \quand
  \tau^{ (i) } (x) \coloneqq \Card \{ 1 \le n \le \tau (x) \colon f^{ n } x \in X_{ i } \}
  ,
\]
and let $F = f^\tau \colon Y \to Y$ be the first return map to $Y$.
For every $I_j \in \cP$ let
\[
  Y_j = 
  \begin{cases}
    I_j \setminus X_{j}, & \text{if } 1 \le j \le d \\
    I_j, &\text{otherwise}. 
  \end{cases}
\]
and for  $i = 1, \ldots, d$ with $i \neq j$ define the sets
\[
  X_{i,n} = f^{-n}_i ( Y ),
  \quand
  Y_{j,i,n} =  f_j^{-1} X_{i,n - 1},
\]
so that $\{ X_{i,n} \}_{ n\ge 1}$ partitions $X_i$, and $\{ Y_{j,i,n} \}_{ i \neq j, n \ge 1}$ partitions $Y_j$. 
Notice that the set $X_{i,n} \subset X_i$ consist of the points that enter $Y$ in exactly $n$ iterations and hence $Y_{j,i,n}  = \{ x \in Y_j \colon \tau^{(i)} (x) = n \}$ and $F|_{Y_{j,i,n}} = f^n \colon Y_{j,i,n} \to Y_i$ is a $C^2$ diffeomorphism.
Hence, letting $\cQ = \{ Y_{j,i,n} \colon i = 1,\ldots, d\;j=1,\ldots, D, \; j \neq i, \; n \ge 1 \}$ we see that $F$ is Markov and, as $F^3 a = Y$ for every $a \in \cQ$, is topologically mixing. Moreover, as $f$ has the Adler property, so does $F$ (see e.g.\ \cite[Lemma 10]{Zwe1998}). As $F$ is uniformly expanding by construction, a standard calculation (see e.g.\ \cite[Proposition 4.3.3]{Aar1997}) 
tells us that the distortion of $F$ is log-Lipschitz. So, we have proved the following Lemma.

\begin{lemma}\label{prop:thaler-gm}
  The first return map $F$ is a topologically mixing Gibbs-Markov map.
\end{lemma}

As $f$ is full branched it topologically mixing, and so it follows then from Lemma~\ref{prop:thaler-gm} and Lemma~\ref{lem:gm-inv-meas} that $f$ preserves a $\sigma$-finite invariant measure absolutely continuous measure whose density $h$ restricted to the set $Y$ is Lipschitz on the elements of $F\cQ = \{ Y_1, \ldots, Y_d \}$. As $\mu$ is the unique absolutely continuous measure for $f$, we know that, after scaling, $d \mu / d m = h$. 

\begin{corollary}
  The map $f$ satisfies \ref{itm:dyn-sep} and \ref{itm:mixing}.
\end{corollary}

\begin{proof}
  As the density $h$ of $\mu$ is bounded on $Y$, we know that $\mu(Y)$ is finite and hence we can assume from hereon that $\mu$ is scaled so that $\mu ( Y )= 1$. Moreover, by construction, the set $Y$ dynamically separates the neighbourhoods $X_1, \ldots, X_d$ of the fixed points $\xi_1, \ldots, \xi_d$ in the sense that \eqref{eq:dyn-sep} holds. This yields \ref{itm:dyn-sep}.

  As the functions $\tau, \tau^{ (1) }, \ldots, \tau^{ (d) }$ are each constant on each of the elements of the Markov partition $\cQ$ we know from Lemma \ref{lem:cf-mix} that \ref{itm:mixing} holds. 
\end{proof}

\begin{lemma}\label{lem:tail}
  \(
  \mu_Y ( \tau^{ ( i ) } = n )
  \sim
  n^{ -2 } \ell_i (n)
  \)
  as $n\to\infty$ for $i=1,\dots,d$,
  where
  \begin{equation}\label{eq:def-ell-tha}
  \ell_i (n) \coloneqq
  u_i \beta_i^\sharp ( 1 / n ) 
  \sum_{1 \le j  \le D \colon j \neq i}
  h ( f_j^{ -1 } \xi_i ) (f_j^{ -1 })' (  \xi_i ),
  \quwhere
  u_i =
  \begin{cases}
    1, & \text{if } \xi_i \in \partial X_i \\ 
    2,  & \text{if } \xi_i \in \operatorname{int} X_i.
  \end{cases}
\end{equation}
\end{lemma}

\begin{proof}
Notice that if the neutral fixed point in the interval $X_i$ lies on the boundary of $X_i$ then the $X_{ i,n }$ each consist of a single interval, otherwise the $X_{i,n}$ are each the union of two disjoint intervals. Moreover, by construction the $X_{i,n}$ accumulate at the corresponding neutral fixed point $\xi_i$. Hence, the sets $Y_{ j,i,n }$ are also the union of at most two intervals which accumulate at the preimage $\zeta_{ j,i }$ of $\xi_{ i }$ in $I_j$. We begin by estimating the size of the $X_{ i,n }$. In order to do so we will make use of the following result based on \cite[Proposition 1]{Hol2005}.

\begin{lemma}\label{lem:hol2005}
  Suppose that $T \colon [0,1] \to \R$ is strictly monotone with $T(0) = 0$ and is given by $T( x ) = x + \gamma (x) x^{ 2 }$ for some function $\gamma$ that is analytic and slowly varying at $0$.  Let $x_0 = T(1)$ and $x_n = T^{ -n } x_{ 0 }$. Then,
  \[
    | x_n - x_{n+1} | \sim n^{ - 2 } \gamma^{ \sharp }( 1 / n ),
  \]
  
  where $\gamma^\sharp$ is the de Bruijn conjugate of $\gamma$.
\end{lemma}

\begin{proof}
  By \cite[Proposition 1]{Hol2005},
  \[
     x_n \sim n^{ - 1 } \gamma^{ \sharp }( 1 / n ).
  \]
  To obtain the asymptotic for the difference we note that
  \[
    T x_{ n + 1 } -  x_{ n + 1 } = \gamma ( x_{ n + 1 } ) x_{ n + 1 }^{ 2 }
    \sim \big[ \gamma \left( (1 / n )  \gamma^{ \sharp } ( 1 / n ) \right) \gamma^{ \sharp } ( 1 / n ) \big]
    n^{ -2 } \gamma^{ \sharp } ( 1 / n )
    \sim 
    n^{ -2 } \gamma^{ \sharp } ( 1 / n )
    ,
  \]
  by the definition of the de Bruijn conjugate.
\end{proof}

Let $i = 1,\ldots, d$.
After an appropriate change of variables we can apply Lemma \ref{lem:hol2005} to each connected component of $X_{i,n}$ to find that length of each component is asymptotic to $ n^{ - 2 } \beta_i^{ \sharp }( 1 / n )$, and hence
\begin{equation}\label{eq:x_n}
  | X_{i,n} | \sim u_i n^{ - 2 } \beta_i^{ \sharp }( 1 / n ).
\end{equation}
Let $1 \le j \le D$ with $j \neq i$. We focus on the case that $X_{i,n}$ has a single connected component, the other case proceeds analogously.
By the mean value theorem there exists a point $y_{ j,i,n } \in Y_{ j,i,n }$ so that $|f'_j y_{ j,i,n }| | Y_{ j,i,n } | = | X_{ i,n } |$. Using the Adler property, a standard calculation (see \cite[Section 4.3]{Aar1997}), shows that $x,y \in I_j$, that $|f'(x) / f'(y)| \ll 1$, and so again by the mean value theorem there is a $z \in I_j$ such that
\[
  \frac{ 1  }{ f' (x) } - \frac{ 1 }{ f'( y ) }
  =
  \frac{ f' (x) - f'( y) }{ f'(x) f'(y) }
  = \frac{ f'' ( z ) }{ f'( z )^2 } \frac{ f'( z ) }{ f' (x) } \frac{ f'( z ) }{ f'( y ) } ( y - x ).
\]
This yields that
\begin{equation}\label{eq:dist} %
  |f'(x)^{ -1 } - f'(y)^{ -1 } | \ll | x - y |,
\end{equation}
where the implied constants are independent of $j,i$ and $n$. Then, using the mean value theorem one more time, we know that there is a $y_j \in I_j$ such that $\frac{ 1 }{ f' ( y_j) } = \frac{ | f ( I_j ) | }{ f'(y_j) } = | I_j |$. Combining this with \eqref{eq:dist} we obtain
\begin{equation}\label{eq:dist2}
    |f'_j y_{ j,i,n }|^{ -1 } \ll | I_j |.
\end{equation}
From \eqref{eq:dist} and \eqref{eq:dist2} we immediately obtain that $|Y_{ j,i,n }| \ll |I_j||X_{ i,n-1 }|$. Also, letting $\eps_n = \sum_{ i = 1 }^{ d } \sum_{ k = n - 1}^{\infty} | X_{i,k} |$,
\[
  | (f' (y_{j,i,n} ) )^{-1}  - (f' (f_{j}^{-1} (\xi_i) ) )^{-1} |
  \ll | y_{j,i,n} - f_{j}^{-1} (\xi_i) |
  \ll \sum_{ k = n }^{ \infty } | Y_{j,i,k} | \ll |I_j| \sum_{ k = n - 1}^{\infty} | X_{i,k} | \ll |I_j| \eps_n,
\]
and so
\[
  |Y_{j,i,n}| = [ f' ( f_j^{-1} (\xi_i)) ]^{ -1 } |X_{i,n-1} | \left[
    1 + O \left(
       |I_j| \eps_n
      \right)
    \right].
\]
As $h$ is Lipschitz on each interval $Y_i$,
\begin{align*}\label{eq:h-calc}
  \big| \mu_Y ( Y_{ j , i, n  } ) - h ( f_{j}^{-1} \xi_i ) | Y_{ j,i,n } | \big|
  &=
    \left|
    \int_{  Y_{ j , i , n }} ( h - h( f_{j}^{-1} \xi_i ) )\,d\Leb
    \right| \\
  &\ll  | Y_{ j , i , n }|\sup_{x \in Y_{j,i,n}}  | x - f_j^{-1} \xi_i| \\ 
  &\ll  | Y_{ j , i , n }| \sum_{ k = n }^\infty | Y_{j,i,k} | \\
  &\ll  | Y_{ j , i , n }| |I_j| \eps_n.
\end{align*}
Thus,
\[
  \mu_Y (Y_{j,i,n}) = h ( f_{j}^{-1} \xi_i ) ' ( f_{j}^{-1} )' (\xi_i) 
  |X_{i,n}| \left[
    1 + O(|I_j| \eps_n)
  \right].
\]
Since $h ( f_{j}^{-1} \xi_i ) ( f_{j}^{-1} )' (\xi_i) \ll | I_j |$ and since $\sum_{ 1 \le j  \le D \colon j \neq i } \eps_n |I_j | \ll \eps_n \to 0$, 
we can sum over the $j \neq i$ to obtain
\[
  \mu_Y ( \tau^{(i)} = n )
  \sim |X_{i,n}|\sum_{1 \le j  \le D \colon j \neq i } h ( f_{j}^{-1} \xi_i ) ( f_{j}^{-1} )' (\xi_i)
  .
\]
Combining this with \eqref{eq:x_n} concludes the proof.
\end{proof}

\begin{proof}
  [Proof of Proposition \ref{prop:thaler}]
  As \ref{itm:f-cons} immediately follows from \eqref{eq:ellnonint}, it remains to check that \ref{itm:fxd-pnts} and \ref{itm:tail} both hold.

  Notice first of all that \ref{itm:tail} follows from Lemma \ref{lem:tail} and the hypotheses \eqref{eq:rho_i} and \eqref{eq:ellnonint}.
 
  Using the definition of the  $\tau^{(i)}$ and \eqref{eq:def-mu}, we see that for any $N$ and any $i = 1 , \ldots, d$,
  \[
    \mu \left( \bigcup_{n = N}^{\infty} X_{i,n} \right)
    \ge \int_{ \{ \tau^{ (i) } \ge N \} } \tau^{ (i) } \,d \mu_Y
    = \int_{N}^\infty \mu_Y (\tau^{ (i) } > x )\,dx
    = \infty,
  \]
  by \eqref{eq:ellnonint}.
  As the sets $X_{i,n}$ (which are a finite union of intervals) are accumulating at the point $\xi_i$ we know that for any $\eps > 0$ there exists an $N$ such that $B_\eps ( \xi_i ) \subset \bigcup_{n = N}^{\infty} X_{i,n}$. Hence $\mu ( B_\eps( \xi_i ) ) = \infty$. This concludes the proof of \ref{itm:fxd-pnts}.
\end{proof}

\subsection{Parabolic Rational Maps}\label{subsec:para-rat-map}
Throughout this section we will let $f$ be a map satisfying \ref{p1} and \ref{p2}.
In this section we will prove Theorem \ref{thm:parabolic}.
We begin with the following lemma which can be deduced from \cite{AarDenUrb1993}.

\begin{lemma}
  The map $g$, after potentially increasing $M$, satisfies \ref{itm:f-cons}--\ref{itm:mixing}.
\end{lemma}
\begin{proof}
  The details of how this result can be deduced from \cite{AarDenUrb1993} are essentially contained in \cite[Example 4.20]{CoaMel2025-sub}. In \cite[Example 4.20]{CoaMel2025-sub} it is shown, after potentially passing to an iterate of $g$,  that there exists a decomposition $X = Y \cup X_1 \cup \cdots \cup X_{ d' }$ so that:
  \begin{itemize}
    \item
    each $X_i$ is an open neighbourhood of the corresponding fixed point;
    
    \item
    the first return map to $Y$ is Gibbs-Markov with respect to a partition whose elements are denoted by $Y_{ n,r }$ where $n \ge 1$ and $r$ is in some finite indexing set $\mathbb{A}$;

    \item
    the return time is equal to $n$ on each $Y_{ n,r }$ and the iterates $g^k Y_{n,r}$ for $k \leq n$ all lie in the same $X_i$ where $i = \phi (r)$ and $\phi \colon \mathbb{A} \to \{1, \ldots, d' \}$ is surjective (and hence $Y$ dynamically separates the fixed points);

    \item there exist constants $c_1,\ldots, c_d$ such that $\mu_Y ( \tau^{ (i) } = n ) = \sum_{ r \colon \phi(r) = i } Y_{ n,r } \sim c_i n^{-2}$ (see \cite[Equation (4.9)]{CoaMel2025-sub}). 
  \end{itemize}
  Combining the above with Lemma \ref{lem:cf-mix} concludes that \ref{itm:f-cons}--\ref{itm:mixing} hold for $g$ and hence concludes the proof.
\end{proof}

Let $\cS' \subset \R^{ d' }$ denote the $d'$ dimensional simplex, let $\vrho' = \frac{ 1 }{ c_1 + \cdots + c_{ d' } } ( c_1, \ldots, c_{ d' })$ and let $\RR_i' \subset \cS'$ be the line that connects the point $\vrho'$ to the $i^{ th }$ vertex. For every $\vp \in \cS'$ we define $\pi ( \vp ) = \sum_{ i = 1 }^{ d' } p_i \delta_{ \xi_i }$.
As $g$ satisfies \ref{itm:f-cons}--\ref{itm:mixing}, and as the slowly varying function that appears in each of the tails $\mu_Y ( \tau^{ (i) } > n )$ is asymptotically constant we know from Corollary \ref{cor:constant} that
\begin{equation}\label{eq:bar_en}
  \{ \text{limit points of } \mu_n \}
  \aseq
  \pi \left( \bigcup_{ i = 1 }^{ d' } \RR_i' \right)
\end{equation}
where $\mu_n$ denote the empirical measures for $g$.

For $i = 1, \ldots, d$ set
\[
  P_i = \frac{ M }{ \text{period of the }i^{ th }\text{ periodic orbit of  }f }
\in \N
  .
\]
Let $\cI = \{(i,j) \colon i = 1, \ldots, d, \; j = 0, \ldots, P_i - 1 \}$ and $\{\xi_{ i,j } \colon (i,j) \in \cI \}$ be a relabelling of the fixed points of $g$ so that each $\xi_{i,j}$ is an element of the $i^{th}$ periodic orbit of $f$ and $f \xi_{i,j} = \xi_{i,j+ 1 \mod P_i }$. In what follows we will also reindex
components of vectors in the simplex
$\cS'$ by $\cI$ in the same way. 

\begin{lemma}\label{lem:periods}
  Let $\vp \in \cS'$ and let $x \in X$ be in the full measure subset for which the limit points $\mu_n \circ f^k$ are contained in $\pi ( \cS' )$ for every $k \ge 0$.
  Suppose that $\mu_{n_m} (x) \to \pi ( \vp )$, then $\mu_{n_m} (f^k x) \to \pi ( \sigma^k \vp )$, where
  $(\sigma^k \vp )_{i,j} = \vp_{i,j - k \mod P_i }$ for all $(i,j) \in \cI$.
\end{lemma}

\begin{proof}
  Let $(i,j) \in \cI$ and let $U_{i,j}$ be an open neighbourhood $\xi_{i,j}$ that contains no other neutral fixed point of $g$ and small enough so that $f^{-1 }U_{i,j} \cap \{ \xi_{ l,q } \colon (l,q ) \in \cI\} = \{ \xi_{i, j - 1 \mod P_i} \}$.
  Then, as $f^{-1 }U_{i,j}$ is an open neighbourhood of $\xi_{i, j - 1 \mod P_i}$ and as $f$ and $g$ commute,
  \[
    \mu_{ n_m } ( f x ) (U_{ i,j }) = 
    \frac{ 1 }{ n_m  }\sum_{ k = 0 }^{ n_m - 1 } 1_{ U_{ i,j } } \circ g^k \circ f ( x )
    =  \frac{ 1 }{ n_m  }\sum_{ k = 0 }^{ n_m - 1 } 1_{ f^{ -1 } U_{ i,j } } \circ g^k ( x )
    = \mu_{n_m} ( x ) ( f^{ -1 }U_{i,j} )
    \to \vp_{i,j - 1 \mod P_i}.
  \]
  Since, $\cL(f x) \subset \pi ( \cS' )$ we know from the above that the only possible limit point of $\mu_{ n_m } ( fx )$ is $\pi ( \sigma \vp )$. Proceeding inductively concludes the proof.
\end{proof}

\begin{proof}[Proof of Theorem \ref{thm:parabolic}]
  For each $n \ge 1$ let $k_n \ge 0$ and $0 \le t_n \le M - 1$ be the integers that satisfy $n = k_n M + t_n$ and notice that $ k_n / n \to 1 / M$ and $t_n / n \to 0$.
  Thus
  \[
    \sum_{ j = 0 }^{ n - 1 } f_*^j \delta_x
    =
    \sum_{ j = 0 }^{ k_n - 1 } g^j_* ( \delta_x + \cdots + \delta_{ f^{ M - 1 } x} ) + \sum_{ j = 0 }^{ t_n - 1 } f_*^{ j+k_nM} \delta_x,
  \]
  and, with $\operatorname{dist}_W$ denoting the Wasserstein distance on $\cM_1(X)$, we have
  \begin{equation*}
    d_W \left( e_{n} (x),  \frac{ 1 }{ M } \left(\mu_{k_n} (x) + \cdots +  \mu_{k_n} ( f^{ M - 1 }x)\right) \right) \ll \frac{ 1 }{ n }.
  \end{equation*}
  Similarly, there exists a subsequence $n_k \sim  k M$ such that
  \begin{equation*}
    d_W \left( e_{n_k} (x),  \frac{ 1 }{ M } \left(\mu_{k} (x) + \cdots +  \mu_{k} ( f^{ M - 1 }x)\right) \right) \ll \frac{ 1 }{ k }.
  \end{equation*}
  Hence, it suffices to describe the limit points of the measure $\bar \mu_k = \frac{ 1 }{ M } \left(\mu_{k}  + \cdots +  \mu_{k} \circ f^{ M - 1 }\right)$.

  Let $x$ be in the full measure set on which \eqref{eq:bar_en} and Lemma \ref{lem:periods} hold.
  By Lemma \ref{lem:periods}, if $\pi\vp$ is a limit point of $\mu_k(x)$ for some $\vp \in \cS$, then $\Pi ( \Psi(\vp) )$ is a limit point of $\bar \mu_k(x)$ where $\Psi(\vp) \in \cS$ is given by $\Psi(\vp)_i = \sum_{ j = 0  }^{ P_i - 1 }  p_{i,j}$. On the other hand, if $\bar \mu_k (x)$ converges to some $\nu$ along a subsequence $k_n$ then, maybe after passing to  another subsequence $k'_n$, the $\mu_{ k'_n}(x)$ must also converge to some limit which necessarily is given by $\pi(\vp)$ for some $\vp \in \cS'$. So, we must have $\nu = \Pi ( \Psi (\vp))$.
  
  The result follows by noticing that $\Psi \left( \bigcup_{ i = 1 }^{ d' } \RR_i' \right) =  \bigcup_{ i = 1 }^{ d } \RR_i$, where
  \begin{equation}\label{eq:def-rho-para}
    \vrho_i = \sum_{ j = 0  }^{ P_i - 1 }  \rho'_{i,j}.
  \end{equation}
\end{proof}

\section*{Acknowledgements}
D.C\ was supported by FAPERJ – Carlos Chagas Filho Foundation for Research Support of the State of Rio de Janeiro, Processes E-26/200.027/2025 and 200.028/2025.

T.S.\ was supported by the European Union’s Horizon Europe research and innovation programme under the Marie Sk{\l}odowska-Curie grant agreement ErgodicHyperbolic -- 101151185.

A visit of D.C.\ to the University of Exeter was supported by a Heilbronn Institute small grant.

The authors would like to thank Zemer Kosloff and Roland Zweim\"uller for useful discussions. 

\printbibliography

\end{document}


%% file: preamble.tex
\usepackage[dvipsnames]{xcolor}
\usepackage{amsmath,amsfonts,mathtools,amssymb}
\usepackage{todonotes}
\usepackage[normalem]{ulem}

\usepackage[url=false,eprint=false,style=alphabetic,isbn=false,sorting=nyt, giveninits=true]{biblatex}

\DeclareFieldFormat{doi}{%
  \href{https://doi.org/#1}{\nolinkurl{#1}}%
}

\renewbibmacro{in:}{%
  \ifentrytype{article}{}{\printtext{\bibstring{in}\intitlepunct}}}

\usepackage{hyperref}
\usepackage{color}
\hypersetup{
 colorlinks=true,
 linkcolor=MidnightBlue,
 citecolor=OliveGreen,
 urlcolor=MidnightBlue
}
\usepackage{comment}

\usepackage{nameref}
\makeatletter
\let\orgdescriptionlabel\descriptionlabel
\renewcommand*{\descriptionlabel}[1]{%
  \let\orglabel\label
  \let\label\@gobble
  \phantomsection
  \edef\@currentlabel{#1}%
  \let\label\orglabel
  \orgdescriptionlabel{#1}%
}

\usepackage{subcaption}
\usepackage{graphicx}

\usepackage{amsthm}

\theoremstyle{definition}
\newtheorem{definition}{Definition}[section]

\theoremstyle{definition}

\theoremstyle{plain}
\newtheorem{lemma}[definition]{Lemma}

\newtheorem{theorem}[definition]{Theorem}
\newtheorem*{theorem*}{Theorem}

\newtheorem{manualtheoreminner}{Theorem}
\newenvironment{manualtheorem}[1]{%
  \renewcommand\themanualtheoreminner{#1}%
  \manualtheoreminner
}{\endmanualtheoreminner}

\newtheorem{proposition}[definition]{Proposition}

\theoremstyle{plain}
\newtheorem{corollary}[definition]{Corollary}

\theoremstyle{remark}
\newtheorem{remark}[definition]{Remark}

\theoremstyle{definition}
\newtheorem{example}[definition]{Example}

\theoremstyle{plain}

\makeatother

\newcommand{\cB}{\mathcal{B}}

\newcommand{\cL}{\mathcal{L}}
\newcommand{\cP}{\mathcal{P}}
\newcommand{\cQ}{\mathcal{Q}}
\newcommand{\cI}{\mathcal{I}}
\newcommand{\cS}{\mathcal{S}}
\newcommand{\cM}{\mathcal{M}}

\newcommand{\eps}{\varepsilon}

\newcommand{\Card}{\operatorname{Card}}
\newcommand{\Leb}{\operatorname{Leb}}

\newcommand{\quand}{\quad\text{and}\quad}
\newcommand{\quwhere}{\quad\text{where}\quad}

\newcommand{\vtau}{\boldsymbol{\tau}}
\newcommand{\vS}{\boldsymbol{S}}

\newcommand{\vrho}{\boldsymbol{\rho}}

\newcommand{\vi}{\boldsymbol{i}}

\newcommand{\vp}{\boldsymbol{p}}

\newcommand{\vv}{\boldsymbol{v}}
\newcommand{\vz}{\boldsymbol{\zeta}}
\newcommand{\vg}{\boldsymbol{\gamma}}
\newcommand{\RR}{\boldsymbol{R}}

\newcommand{\R}{\mathbb{R}}
\newcommand{\C}{\mathbb{C}}
\newcommand{\N}{\mathbb{N}}
\renewcommand{\P}{\mathbb{P}}

\newcommand{\E}{\mathbb{E}}

\newcommand{\aseq}{\overset{as}{=}}

\newcommand{\dist}{\operatorname{dist}}
\newcommand{\NN}{\mathfrak{N}}

\definecolor{eggshell}{rgb}{0.94, 0.92, 0.84}